\documentclass{amsart}
\usepackage{amssymb,latexsym}
\input epsf
\usepackage{graphicx}
\newtheorem{theorem}{Theorem}[section]
\newtheorem{proposition}[theorem]{Proposition}
\newtheorem{corollary}[theorem]{Corollary}
\newtheorem{definition}[theorem]{Definition}
\newtheorem{remark}[theorem]{Remark}
\newtheorem{lemma}[theorem]{Lemma}

\newcommand{\nc}{{\rm {\bf L}}}

\newcommand{\fF}{{\mathcal F}}
\newcommand{\mM}{{\mathcal M}}
\newcommand{\rR}{{\mathcal R}}
\newcommand{\RR}{{\mathbb R}}

\newcommand{\sm}{{\setminus}}
\begin{document}
\title[Generalized associahedra and noncrossing partitions]
{$h$-vectors of generalized associahedra and noncrossing partitions}

\author[Athanasiadis]{Christos~A.~Athanasiadis}
\address{Department of Mathematics
(Division of Algebra-Geometry)\\
University of Athens\\
Panepistimioupolis\\
15784 Athens, Greece}
\email{caath@math.uoa.gr}

\author[Brady]{Thomas~Brady}
\address{School of Mathematical Sciences\\
Dublin City University\\
Glasnevin, Dublin 9\\
Ireland}
\email{tom.brady@dcu.ie}

\author[McCammond]{Jon~McCammond}
\address{Mathematics Department\\
University of California, Santa Barbara\\
Santa Barbara, CA 93106\\
U.S.A} \email{jon.mccammond@math.ucsb.edu}

\author[Watt]{Colum~Watt}
\address{School of Mathematical Sciences\\
Dublin Institute of Technology\\
Dublin 8\\
Ireland}
\email{colum.watt@dit.ie}

%%\subjclass{Primary ; Secondary .}
\date{February 6, 2006}
%%\\Supported in part by the American Institute of Mathematics (AIM)
%%and the NSF}
%
\begin{abstract}
A case-free proof is given that the entries of the $h$-vector of
the cluster complex $\Delta (\Phi)$, associated by S. Fomin and A.
Zelevinsky to a finite root system $\Phi$, count elements of the
lattice $\nc$ of noncrossing partitions of corresponding type by
rank. Similar interpretations for the $h$-vector of the positive
part of $\Delta (\Phi)$ are provided. The proof utilizes the
appearance of the complex $\Delta (\Phi)$ in the context of the
lattice $\nc$, in recent work of two of the authors, as well as an
explicit shelling of $\Delta (\Phi)$.
\end{abstract}

\maketitle

\section{Introduction}
\label{intro}

Let $\Phi$ be a finite root system of rank $n$ with corresponding
finite real reflection group $W$. The cluster complex $\Delta
(\Phi)$ was introduced by S. Fomin and A. Zelevinsky in the
context of algebraic $Y$-systems \cite{FZ2} and cluster algebras
\cite{FZ1, FZ3}. It is a pure $(n-1)$-dimensional simplicial
complex which is homeomorphic to a sphere \cite{FZ2}. When $\Phi$
is crystallographic there is a cluster algebra associated to
$\Phi$ and $\Delta (\Phi)$ was realized explicitly in \cite{CFZ}
as the boundary complex of a simplicial convex polytope, known as
a (simplicial) \emph{generalized associahedron}; see \cite{FR} for
an overview of cluster complexes and generalized associahedra.

The combinatorics of $\Delta (\Phi)$ encodes the exchange of clusters in
the corresponding cluster algebra of finite type. The remarkable enumerative
properties of $\Delta (\Phi)$ have provided combinatorialists with many
interesting puzzles and have suggested connections with a number of seemingly
unrelated objects appearing in various areas of mathematics; see \cite{Ch,
FZ2} and \cite[Lecture 5]{FR}. One such object is the lattice of noncrossing
partitions associated to $W$ (see \cite{ABW, Be, BW1, BW, McC} and Section
\ref{nc}). This is a self-dual graded poset of rank $n$, denoted $\nc_W$,
which is a lattice of considerable interest in the topology of
finite-type Artin groups. It encodes the combinatorics of
reduced decompositions of Coxeter elements of $W$ with respect to the
generating set of all reflections in $W$. At the core of the connection between
$\Delta (\Phi)$ and $\nc_W$ lies the fact that both the number of facets of
$\Delta (\Phi)$ (clusters) and the number of elements of $\nc_W$ are given
by the expression
\begin{equation}
N (\Phi) = \prod_{i=1}^n \frac{e_i + h + 1}{e_i + 1},
\label{prod}
\end{equation}
known as the \emph{Catalan number} associated to $W$, where $h$ is the
Coxeter number and $e_1,\dots,e_n$ are the exponents of $W$. Recall (see
Section \ref{complexes} for definitions) that the $h$-vector of a
simplicial complex $\Delta$ is a fundamental enumerative invariant which
refines the number of facets of $\Delta$. The following statement has been
verified by use of the classification of finite root systems and case by
case computations due to various authors; see parts (i) and (ii)
of \cite[Theorem 5.9]{FR}.

\begin{theorem}
Let $\Phi$ be a finite root system of rank $n$ with corresponding reflection
group $W$. The entry $h_i (\Delta (\Phi))$ of the $h$-vector of $\Delta
(\Phi)$ is equal to the number of elements of $\nc_W$ of rank $i$ for all $0
\le i \le n$.

In particular, the number of facets of $\Delta (\Phi)$ is equal to the total
number of elements of $\nc_W$.
\label{thm0}
\end{theorem}

It has been one of the challenging problems on cluster
combinatorics to find a conceptual proof of the previous theorem.
That such a proof may be possible is suggested by the recent work
\cite{BW} of two of the authors.  There, a case-free proof of
the lattice property of $\nc_W$ is given together with a new
characterization of the cluster complex. It is the main objective
of this paper to build upon the constructions of \cite{BW} and
give a case-free proof of Theorem \ref{thm0}. More specifically we
will describe a shelling of $\Delta (\Phi)$ and a bijection $\phi$
from the set of facets of $\Delta (\Phi)$ to $\nc_W$ such that,
for any facet $F$ of $\Delta (\Phi)$, the number of vertices of
the restriction face of $F$ with respect to this shelling is equal
to the corank of $\phi (F)$ in $\nc_W$. Our shelling of $\Delta
(\Phi)$ is a variation of the lexicographic ordering on the facets
with respect to the total ordering of the vertices of $\Delta
(\Phi)$ considered in \cite{BW}. The fact that the map $\phi$ (see
Definition \ref{def:phi}) is indeed a bijection is the hardest
part of the proof and requires a thorough analysis of certain
properties of the subcomplexes $X(w)$ of $\Delta (\Phi)$ defined
and studied in \cite{BW}. For another conceptual approach to
Theorem \ref{thm0}, see \cite{Re, RS}.

The entries of the $h$-vector
of $\Delta (\Phi)$ are known as the \emph{Narayana numbers} associated to
$W$. They are the dimensions of
the real cohomology groups of the complex projective toric variety associated
to the corresponding generalized associahedron and they admit
various interesting combinatorial interpretations (one of which appears
in Theorem \ref{thm0}) and generalizations; see \cite[Section 5.2]{FR}
and the references given there.

\medskip
The general layout of this paper is as follows. Basic definitions and
background related to noncrossing partition lattices and cluster complexes
are collected in Section \ref{pre}. In particular, the realization of $\Delta
(\Phi)$ given in \cite{BW} and other necessary material from that paper
are reviewed. Some results from \cite{BW} are also generalized and some new
technical results are given. The description \cite[Section 6]{BW} of the
lexicographically first facet of the complex $X(w)$, with respect to the
vertex ordering of \cite{BW}, is recalled in Section \ref{lex} and an
analogous description for the lexicographically last facet is provided.
These constructions are used in establishing some key properties (Theorem
\ref{thm:lexlast} and Proposition \ref{prop:per}) of the map $\phi$.
This map is defined in Section \ref{types} after some further definitions and results
related to the structure of faces of $\Delta(\Phi)$ and the complexes $X(w)$
are given in Section \ref{sec:sub}.
The fact that $\phi$ is a bijection is
proved in Section \ref{bij}. The proof of Theorem \ref{thm0} is completed in
Section \ref{shell} by combining bijectivity of $\phi$ with the explicit
shelling of $\Delta (\Phi)$ mentioned above. Some combinatorial
interpretations of the entries of the $h$-vector of the positive part of
$\Delta (\Phi)$ are also derived. General background on root systems, cluster
complexes and noncrossing partition lattices can be found in \cite{FR} and
references therein.

\section{Preliminaries}
\label{pre}

Throughout this paper $W$ is a finite real reflection group of rank $n$,
generated by orthogonal reflections in $\RR^n$. Unless otherwise stated  we
assume that $W$ is irreducible. The reflection
in the linear hyperplane in $\RR^n$ orthogonal to a nonzero vector $\alpha$
is written as $R(\alpha)$.  The root system for $W$ which we use is
denoted by $\Phi$ and consists of the pairs
$\{\alpha, -\alpha\}$ of normals to the reflecting hyperplanes.  We will further
assume that these normals have \emph{unit} length.

\subsection{The reflection length order and the lattice $\nc_W$.}
\label{nc}

To any $w \in W$ we can associate two linear subspaces of $\RR^n$: the fixed
space $\fF(w)$ of $w$ and the orthogonal complement $\mM(w)$ of $\fF(w)$ in
$\RR^n$.  The subspace  $\mM(w)$ is called the \emph{moved space} of
$w$ in \cite{BW}.  The dimension $\ell(w)$
of $\mM(w)$ is equal to the smallest integer $k$ such that $w$ can be
written as a product of $k$ reflections in $W$ \cite[Proposition
2.2]{BW1}.   Note that $\ell$ is \emph{not} the usual length function
associated to a simple system for $W$.  For $a, b \in W$ we let
\[ a \preceq b \ \ \ \text{if and only if} \ \ \ \ell (a) + \ell (a^{-1}
b) = \ell(b). \]
The relation $\preceq$ turns $W$ into a graded partially ordered set of rank
$n$, which has the identity $I$ as its unique minimal element and rank function
$\ell$. The following lemma collects some useful properties of this partial
order.

\begin{lemma}
Let $a, b, c \in W$.
\begin{enumerate}
\itemsep=0pt
\item[(i)] If $a \preceq b$ then $\mM(a) \subseteq \mM(b)$.
\item[(ii)] If $\mM(a) \subseteq \mM(b)$ and $a$ is a reflection then $a
\preceq b$.
\item[(iii)] If $\mM(a) \subseteq \mM(b)$ and $a, b \preceq w$ for some
$w \in W$ then $a \preceq b$.
\item[(iv)] If $a, b \preceq c \preceq w$ for some $w \in W$ and $ab \preceq
w$ then $ab \preceq c$.
\end{enumerate}
\label{lem:basic1}
\end{lemma}
\begin{proof}
For parts (i), (ii) and (iii) see \cite[Section 2]{BW1}. Part (iv) follows
from part (iii) since from the assumptions we have $\mM(a), \mM(b) \subseteq
\mM(c)$ and hence $\mM(ab) \subseteq \mM(a)+\mM(b) \subseteq \mM(c)$.
\end{proof}

If $\gamma$ is a Coxeter element of $W$ then the interval $[1, \gamma]$ in
the poset $(W, \preceq)$, denoted $\nc_W (\gamma)$ or $\nc (\gamma)$, is
a locally self-dual, graded lattice of rank $n$; see \cite{ABW, Be, BW1, BW}
for further information and interesting properties. The isomorphism type
$\nc_W$ of $\nc_W (\gamma)$ is independent of $\gamma$ and called the
\emph{noncrossing partition lattice} associated to $W$.

\subsection{Coxeter elements and the operator $\mu$.}
\label{gamma}

As in \cite{Be, BW} our standard choice of $\gamma$, which we fix once and
for all, is that of a bipartite Coxeter element. Thus we fix an
ordered simple system $\Pi = \{\alpha_1, \alpha_2,\dots,\alpha_n \}$ for
$\Phi$ with the property that there is a $1 \le s \le n$ for which
$\Pi_1 = \{\alpha_1,\dots,\alpha_s\}$ and $\Pi_2
= \{\alpha_{s+1},\dots,\alpha_n\}$ are orthonormal sets and we set
\[ \gamma = R(\alpha_1) R(\alpha_2) \cdots R(\alpha_n). \]
Following \cite{BW} we write
\[ \mu(x) = 2 (I - \gamma)^{-1} (x) \]
for $x \in \RR^n$ and recall from \cite[Corollary 3.3]{BW} that for $\tau
\in \Phi$, $\mu (\tau)$ is the unique vector in the one dimensional space
$\fF (R(\tau) \gamma)$ satisfying $\mu (\tau) \cdot \tau = 1$. The following
lemma is implicit in \cite{BW}.
\begin{lemma}
For nonparallel roots $\sigma, \tau$ we have $R(\sigma) R(\tau) \preceq
\gamma$ if and only if $\mu(\sigma) \cdot \tau = 0$.
\label{lem:basic}
\end{lemma}
\begin{proof}
The statement $R(\sigma) R(\tau) \preceq \gamma$ is equivalent to $R(\tau)
\preceq R(\sigma) \gamma$. By parts (i) and (ii) of Lemma \ref{lem:basic1},
this in turn holds if and only if $\tau$ lies in the $(n-1)$-dimensional
space $\mM (R(\sigma) \gamma)$ or, in other words, if and only if
$\mu(\sigma) \cdot \tau = 0$.
\end{proof}
Lemma \ref{lem:sigma-mu} and part (i) of Lemma \ref{lem:tau-rho} below
generalize \cite[Theorem 3.2]{BW} and \cite[Theorem 3.7 (a)]{BW},
respectively.
\begin{lemma}
If $\tau \in \Phi$ and $R(\tau) \preceq w \preceq \gamma$ then $w(\mu(\tau))
= \mu(\tau)-2\tau$.
\label{lem:sigma-mu}
\end{lemma}
\begin{proof}
Observe that $\mu(\tau) \in \fF (R(\tau)\gamma) \subseteq
\fF (w^{-1}\gamma)$, since $w^{-1} \gamma \preceq R(\tau) \gamma$. Hence
$w(\mu(\tau)) = \gamma(\mu(\tau)) = \mu(\tau)-2\tau$.
\end{proof}
\begin{lemma}
If $\sigma, \tau \in \Phi$ and $R(\sigma), R(\tau) \preceq w \preceq \gamma$
then
\begin{enumerate}
\itemsep=0pt
\item[(i)] $w(\mu(\sigma))\cdot \tau = -\mu(\tau) \cdot \sigma$,
\item[(ii)] $w(\mu(\sigma))\cdot \tau = \mu(w(\sigma))\cdot \tau$,
\item[(iii)] $\mu(w(\sigma)) \cdot w(\tau) = \mu(\sigma) \cdot \tau$.
\end{enumerate}
\label{lem:tau-rho}
\end{lemma}
\begin{proof}
For part (i) use Lemma \ref{lem:sigma-mu} to write $\tau = (1/2)
\,[\mu(\tau) - w(\mu(\tau))]$, and similarly for $\sigma$, so that
\begin{eqnarray*}
w(\mu(\sigma))\cdot \tau &=& w(\mu(\sigma))\cdot
(1/2) \, [\mu(\tau)-w (\mu(\tau))] \\
&=& (1/2) \, [w(\mu(\sigma))\cdot \mu(\tau)-w(\mu(\sigma))
\cdot w (\mu(\tau))] \\
&=& (1/2) \, [w(\mu(\sigma))\cdot \mu(\tau)-\mu(\sigma)\cdot \mu(\tau)] \\
&=& (1/2) \, \mu(\tau) \cdot [w(\mu(\sigma))-\mu(\sigma)] \\
&=& -\mu(\tau) \cdot \sigma.
\end{eqnarray*}
For part (ii) observe that $R(w(\sigma)) \preceq w$ and apply part (i) twice
to get
\[\mu(w(\sigma)) \cdot \tau = -w(\mu(\tau)) \cdot w(\sigma)
= - \mu(\tau) \cdot \sigma = w(\mu(\sigma))\cdot \tau. \]
Similarly part (iii) follows by applying part (ii) with $\tau$ replaced by $w(\tau)$.
\end{proof}
%
%\textbf{Aside:  }  If $\sigma = \gamma$  and $\rho$, $\tau$ are roots then
%\[\mu(\rho) \cdot \tau + \mu(\tau) \cdot \rho = 2 \rho \cdot \tau\]
%\vskip .2cm

\subsection{The ordering of roots.}
\label{roots}

Let us denote by $\Phi^+$ the positive system of $\Phi$ corresponding to
the simple system $\Pi$ and by $N$ the cardinality of $\Phi^+$ and let us
use the notation $\Phi_{\ge -1}$ of \cite{FZ2} for the set $\Phi^+ \cup
(-\Pi)$ of almost positive roots in $\Phi$.

As in \cite{BW} we set $\rho_i = R(\alpha_1) R(\alpha_2) \cdots
R(\alpha_{i-1}) (\alpha_i)$ for $i \ge 1$ where the $\alpha_i$ are indexed
cyclically modulo $n$ (so that $\rho_1 = \alpha_1$) and $\rho_{-i} =
\rho_{2N-i}$ for $i \ge 0$ and recall that
\[ \begin{tabular}{ll}
$\{ \rho_1, \rho_2,\dots,\rho_N \} = \Phi^+$, \\
$\{ \rho_{N+i}: 1 \le i \le s\} = \{-\rho_1,\dots,-\rho_s \} =
- \Pi_1$, \\
$ \{\rho_{-i}: 0 \le i < n-s \} = \{-\rho_{N-i}: 0 \le i < n-s\} =
- \Pi_2$.
\end{tabular} \]
Thus, as in \cite{BW}, we can consider the total order $<$ of the set
$\Phi_{\ge -1}$ defined by
\begin{equation}
\rho_{-n+s+1} < \cdots < \rho_0 < \rho_1 < \cdots < \rho_{N+s}.
\label{total}
\end{equation}
Unless otherwise stated, whenever we talk
about the lexicographic order on a collection of subsets of
$\Phi_{\ge -1}$ it should be understood that this is defined with
respect to the total order (\ref{total}). The set
$\{\rho_{N-n+1},\dots,\rho_N\}$ of the last $n$ positive roots in
this order will be denoted by $\Omega$. The following three lemmas from
\cite{BW} are listed here for easy reference.
\begin{lemma} {\rm (\cite[Lemma 3.9]{BW})}
For a set $\{\tau_1, \tau_2,\dots,\tau_k\}$ of
positive roots satisfying $\tau_1 < \tau_2 < \dots < \tau_k$ the following are
equivalent:
\begin{enumerate}
\itemsep=0pt
\item[(a)] $R(\tau_k) \cdots R(\tau_1)$ is an element of $\nc (\gamma)$
of length $k$,
\item[(b)] $\mu(\tau_j) \cdot \tau_i = 0$ for $1 \le i < j \le k$. \qed
\end{enumerate}
\label{lem:BW}
\end{lemma}
\begin{lemma} {\rm (\cite[Theorem 3.7 (b) and (d)]{BW})}
For $1 \le i < j \le N$ we have
\begin{enumerate}
\itemsep=0pt
\item[(i)] $\mu(\rho_i) \cdot \rho_j \ge 0$ and
\item[(ii)] $\mu(\rho_j) \cdot \rho_i \le 0$. \qed
\end{enumerate}
\label{lem:BW2}
\end{lemma}
\begin{lemma} {\rm (\cite[Lemma 5.6]{BW})}
Let $\tau, \rho$ be distinct positive roots with $R(\tau) R(\rho)
\preceq \gamma$.
\begin{enumerate}
\itemsep=0pt
\item[(i)] If $\tau < \rho$ then $\tau \cdot \rho \le 0$.
\item[(ii)] If $\tau > \rho$ then $\tau \cdot \rho \ge 0$. \qed
\end{enumerate}
\label{lem:BW3}
\end{lemma}
\begin{remark}
{\rm  Recall from \cite[Section 3]{BW} that $\gamma^{-1} = R(\rho_i) R(\rho_{i+1})
\cdots R(\rho_{i+n-1})$ and that $\gamma (\rho_i) = \rho_{i+n}$ for $i \ge
1$. In particular, $\Omega$ consists of those positive roots $\rho$ such
that $\gamma (\rho)$ is a negative root. }
\label{rem:lastn}
\end{remark}

The following technical result
gives some information about the action of an arbitrary $w \in \nc (\gamma)$ on $\Phi$.
\begin{corollary}
\begin{enumerate}
\itemsep=0pt
\item[(i)] If $R(\rho) \preceq w \preceq \gamma$ and $\rho$ and $w (\rho)$
are positive roots then $\rho < w(\rho)$.
\item[(ii)] If $R(\rho_i) \preceq w \preceq \gamma$ for some $i \ge 1$ then
$w (\rho_i) \not \in \{\rho_{i+1}, \rho_{i+2},\dots,\rho_{i+n-1}\}$.
\item[(iii)] If $R(\rho) \preceq w \preceq \gamma$ and $\rho \in \Omega$
then $w (\rho)$ is a negative root.
\end{enumerate}
\label{lem:nsteps}
\end{corollary}
\begin{proof}
(i) Using Lemma \ref{lem:tau-rho} we compute $\mu (w(\rho)) \cdot \rho =
- \mu (\rho) \cdot \rho = -1$, so the result follows from part (i) of Lemma
\ref{lem:BW2}.

\noindent
(ii) Suppose that $w (\rho_i) = \rho_j$ with $i < j \le i+n-1$. Then
$R(\rho_j) R(\rho_i) \preceq \gamma$ by Remark~\ref{rem:lastn} and
hence $\mu(\rho_j) \cdot \rho_i = 0$ by Lemma \ref{lem:basic}.
However $\mu(\rho_j) \cdot \rho_i = \mu
(w(\rho_i)) \cdot \rho_i = -1$ as in part (i), which gives a contradiction.

\noindent
(iii) This follows from parts (i) and (ii).
\end{proof}

\subsection{Peripheral elements.}
\label{peripheral}
\begin{definition}
An element $w \in \nc (\gamma)$ is called peripheral if $w \preceq \gamma'$ where
\[ \gamma' \in \{R(\alpha_1)\gamma,\dots,R(\alpha_s)\gamma,
\gamma R(\alpha_{s+1}),\dots,\gamma R(\alpha_n)\}.  \]
Otherwise $w$ is called non-peripheral.
\end{definition}

The following characterizations of peripheral elements of $\nc (\gamma)$ will be useful.
Recall that a standard parabolic subgroup of $W$ is a subgroup generated by a subset
of $\{ R(\alpha_1), \dots , R(\alpha_n)\}$.
\begin{proposition}  If $w \in \nc (\gamma)$ then the following are equivalent:
\begin{enumerate}
\itemsep=0pt \item[(i)] $w$ is peripheral,
\item[(ii)] $w$ lies in
a proper standard parabolic subgroup,
\item[(iii)] at least one of
the roots in $\Omega$ belongs to $\mM (w^{-1} \gamma)$.
\end{enumerate}
\label{prop:per-char}
\end{proposition}
\begin{proof}
(i) $\Rightarrow$ (ii):  If $w$ is peripheral and $R \preceq W$ is a reflection,
then $R$ lies in the standard parabolic subgroup with simple system $\Pi \setminus
\{\alpha_i\}$ for some $1 \le i \le n$.  Thus $w$, which is a product
of such reflections, must belong to the same parabolic subgroup.

\vskip .2cm (ii) $\Rightarrow$ (i): If $w$ belongs to a proper
standard parabolic subgroup then it belongs to one with simple
system $\Pi \setminus \{\alpha_i\}$ for some $1 \le i \le n$ and
$\mM(w) \subseteq \mM(\gamma')$, where
\[ \gamma' \in \{R(\alpha_1)\gamma,\dots,R(\alpha_s)\gamma,
\gamma R(\alpha_{s+1}),\dots,\gamma R(\alpha_n)\}.  \] By part
(iii) of Lemma~\ref{lem:basic1} we have $w \preceq \gamma'$.

\vskip .2cm (i) $\Leftrightarrow$ (iii): Since $R(\alpha_i)\gamma
= \gamma R(\gamma^{-1}(\alpha_i))$ for $1 \le i \le s$, we have
that $w$ is peripheral if and only if $w \preceq \gamma'$ where $\gamma'$ is one of
the elements
\[ \gamma R(\gamma^{-1}(\alpha_1)),\dots,\gamma R(\gamma^{-1}(\alpha_s)),
\gamma R(\alpha_{s+1}),\dots,\gamma R(\alpha_n). \]
This happens if and only if $R \preceq w^{-1} \gamma$ where $R$ is one of
the reflections
\[R(\gamma^{-1}(\alpha_1)),\dots,R(\gamma^{-1}(\alpha_s)),
R(\alpha_{s+1}),\dots,R(\alpha_n). \]
By the definition of the total order (\ref{total}), these are the
reflections defined by the last $n$ positive roots, i.e., defined by the
elements of $\Omega$.
\end{proof}
\subsection{$h$-vectors and shellings.}
\label{complexes}

We will use the notion of a spherical simplicial complex in $\RR^n$.
The faces of such a complex are formed by intersecting the unit sphere $S^{n-1}$
in $\RR^n$ with simplicial cones, each pointed at the origin.  Thus a polyhedral
fan in $\RR^n$ is formed. If $C$ is such a cone of dimension
$n$ and $H$ runs through the supporting hyperplanes of the facets of $C$,
we will refer to the intersections $S^{n-1} \cap H$ as the \emph{walls}
of the spherical simplex $S^{n-1} \cap C$.

Given a finite (abstract, geometric or spherical) simplicial complex
$\Delta$ of dimension $d-1$, let $f_i$ denote the number of
$i$-dimensional faces of $\Delta$. The $h$-\emph{vector} of $\Delta$
is the sequence $h(\Delta) = (h_0, h_1,\dots,h_d)$ defined by
\begin{equation}
\sum_{i=0}^d \, f_{i-1} (x-1)^{d-i} \ = \ \sum_{i=0}^d \, h_i x^{d-i}
\label{h-vec}
\end{equation}
where $f_{-1} = 1$ unless $\Delta$ is empty. The complex $\Delta$ is
\emph{pure} if all its facets (faces which are maximal with respect to
inclusion) have dimension $d-1$. Given a total ordering $F_1,
F_2,\dots,F_m$ of the facets of a pure simplicial complex $\Delta$ of
dimension $d-1$ and $1 \le j \le m$ we denote by $\rR(F_j)$ the set of
vertices $x$ of $F_j$ for which the codimension one face of $F_j$ not
containing $x$ is contained in at least one of the facets $F_1,
F_2,\dots,F_{j-1}$. Such an ordering $F_1, F_2,\dots,F_m$ is called a
\emph{shelling} of $\Delta$ if there are no indices $1 \le i < j \le m$
for which $\rR(F_j)$ is contained in the vertex set of $F_i$. In that case
$\rR(F_j)$ is called the \emph{restriction set} of $F_j$ with respect to
this shelling and the entries $h_i = h_i (\Delta)$ of $h(\Delta)$ are
nonnegative integers given by
\begin{equation}
h_i = \# \ \{1 \le j \le m: \ \#\rR(F_j) = i \}, \ \ \ 0 \le i \le d,
\label{h-shell}
\end{equation}
where the cardinality of a finite set $S$ is denoted by $\# S$.
More information and references on shellability of simplicial
complexes can be found in \cite[Section 11]{Bj}.

\subsection{Cluster complexes and subcomplexes.}
\label{clusters}

Let $\Delta (\gamma)$ denote the collection of spherical simplices in
$\RR^n$ on the vertex set $\Phi_{\ge -1}$ defined by declaring a
subset $\{\tau_1, \tau_2,\dots,\tau_k\}$ of
$\Phi_{\ge -1}$ satisfying $\tau_1 < \tau_2 < \dots < \tau_k$ to be the vertex set
of a simplex in $\Delta (\gamma)$ if and only if
\[ R(\tau_k) R(\tau_{k-1}) \cdots R(\tau_1) \]
is an element of $\nc (\gamma)$ of rank $k$.   Thus $\Delta (\gamma)$
coincides with the spherical simplicial complex $EX(\gamma)$ of \cite{BW}.
A simplex in $\Delta
(\gamma)$ is said to be \emph{positive} if its vertices are
positive roots. For $w \in \nc(\gamma)$ let $X(w)$ denote the
subcollection of $\Delta (\gamma)$ consisting of those simplices
with vertex set contained in $\mM (w) \cap \Phi^+$; in particular,
$X(\gamma)$ is the subcollection of positive simplices of $\Delta
(\gamma)$. The set of vertices (zero dimensional simplices) of $X(w)$
is the positive system induced
by $\Phi^+$ on the root system $\Phi (w) = \Phi \cap \mM(w)$ and
is denoted by $\Phi^+ (w)$. By parts (i) and (ii) of Lemma
\ref{lem:basic}, this set coincides with the set of positive roots
$\tau$ satisfying $R(\tau) \preceq w$. We note that $R(\tau) \rho
\in \Phi (w)$ whenever $\tau, \rho \in \Phi (w)$. Finally, let
$\Delta_+ (\Phi)$ denote the induced subcomplex of $\Delta (\Phi)$
on the vertex set $\Phi^+$, referred to as the \emph{positive
part} of $\Delta (\Phi)$. The following theorem will be crucial in
relating the combinatorics of $\Delta (\Phi)$ to that of $\nc_W$.
\begin{theorem} {\rm (\cite{BW})}
\begin{enumerate}
\itemsep=0pt
\item[(i)] The collection $\Delta (\gamma)$ is a spherical simplicial
complex of dimension $n-1$ which is a realization of the cluster
complex $\Delta (\Phi)$.
\item[(ii)] If $w \in \nc(\gamma)$ then the collection $X(w)$ is a
spherical simplicial complex of dimension $\ell(w)-1$.
\end{enumerate}
In particular, $X(\gamma)$ is a spherical simplicial complex of
dimension $n-1$ which is a realization of $\Delta_+ (\Phi)$. \qed
\label{thm:BW}
\end{theorem}

\bigskip
\begin{figure}
\epsfysize = 5.4 in \centerline{\epsffile{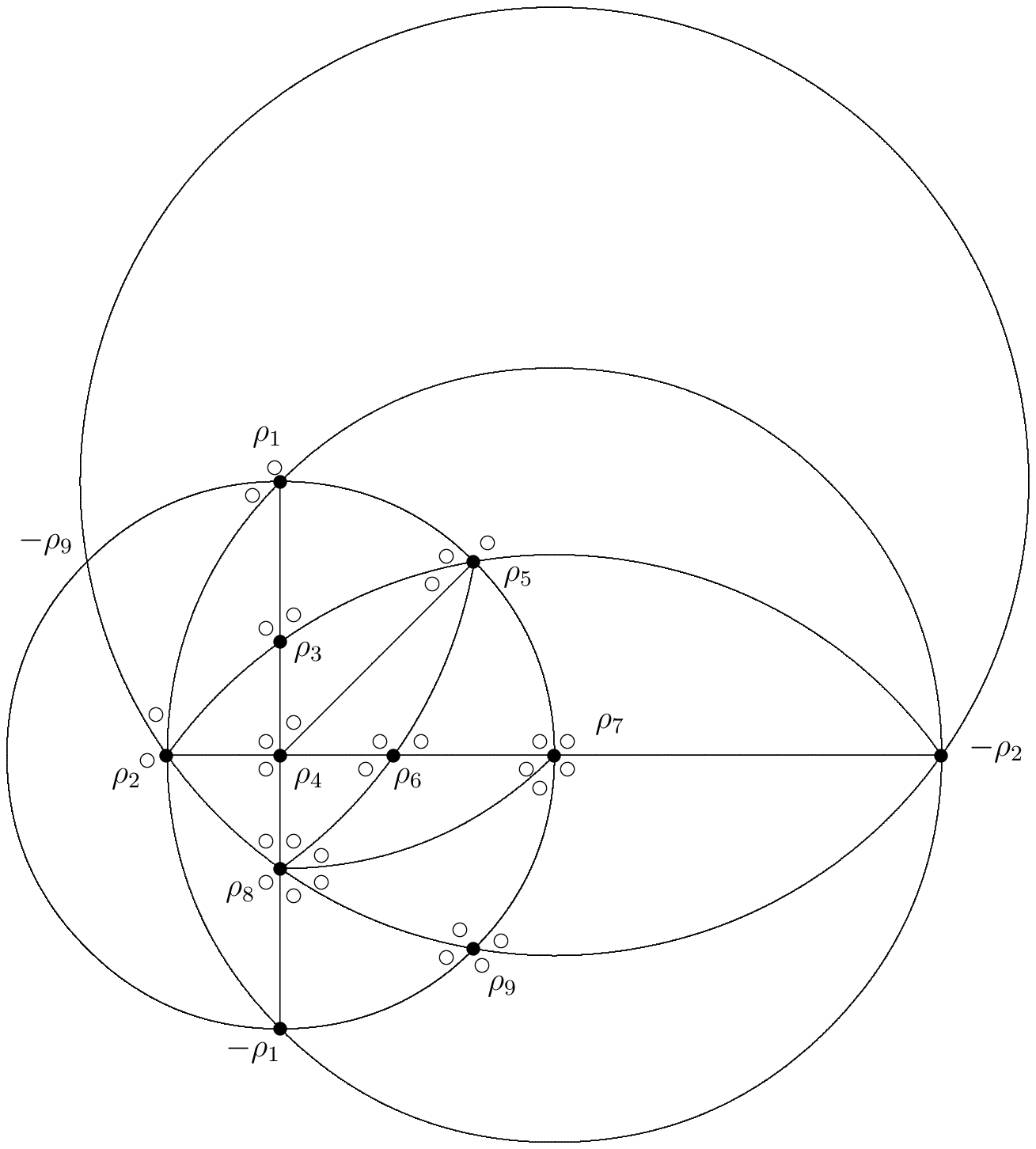}} \caption{}
\label{f:cyclohedron}
\end{figure}

It follows that $\Delta (\gamma)$ and $X(\gamma)$ are pure
simplicial complexes which are homeomorphic to a sphere and a
ball, respectively, of dimension $n-1$.  More generally, for $w
\in \nc(\gamma)$ the complex $X(w)$ is a triangulation of a
spherical simplex of dimension $\ell(w)-1$ \cite[Corollary
7.7]{BW} and hence $X(w)$ is a pure simplicial complex which is
homeomorphic to a ball of dimension $\ell(w)-1$. These facts will
be used repeatedly in Sections \ref{sec:sub} and \ref{shell}.
Figure~\ref{f:cyclohedron} shows $\Delta(\gamma)$ in the case
where $W$ is the group $C_3$ (or $B_3$) of symmetries of the cube,
with the roots ordered as in (\ref{total}), so that $ \rho_0 =
-\rho_9$, while $\rho_{10} = -\rho_1$ and $\rho_{11} = -\rho_2$.

\section{Lexicographically first and last facets}
\label{lex}

In this section we fix $w \in \nc (\gamma)$ of length $k$ and describe
the first and last facets of the complex $X(w)$ in the lexicographic
order. We let $\Pi (w) = \{\delta_1,\dots,\delta_k\}$ be the simple
system for $\Phi^+ (w)$, as
in \cite[Section 5]{BW}, with  $\delta_1 < \delta_2 < \dots < \delta_k$ and
$w = R(\delta_1) \cdots R(\delta_k)$.

\subsection{The lexicographically first facet}
\label{first}

We first recall the description of the first facet of $X(w)$ in the
lexicographic order from \cite[Section 6]{BW}. We provide proofs of
statements, somewhat simplified from those of \cite[Section 6]{BW}, to
make this paper more self-contained and since proofs of the
corresponding statements about the lexicographically last facet will be
similar. We define the roots
\begin{equation}
\epsilon_i = R(\delta_1) \cdots R(\delta_{i-1}) \, \delta_i
\label{def:epsilon}
\end{equation}
for $1 \le i \le k$, so that $\epsilon_1 = \delta_1$. Since $\epsilon_i
- \delta_i$ is a nonnegative linear combination of
$\{\delta_1,\dots,\delta_{i-1}\}$, the set $\{\epsilon_1,\dots,\epsilon_k\}$
is a linearly independent subset of $\Phi^+ (w)$.
\begin{lemma}
We have
\[ \mu(\epsilon_i) \cdot \delta_j =
\begin{cases} 1 & \text{if \ $i=j$} \\
0 & \text{if \ $i\neq j$}. \end{cases} \]
In particular $\mu(\epsilon_i) \cdot \tau \ge 0$ for all $i$ and for all $\tau \in
\Phi^+ (w)$.
\label{lem:epsilon-delta}
\end{lemma}
\begin{proof}
From (\ref{def:epsilon}) we deduce that
\begin{equation}
w = R(\epsilon_i) R(\delta_1) \cdots R(\delta_{i-1}) R(\delta_{i+1})
\cdots R(\delta_k).
\label{eq:epsilon2}
\end{equation}
Since $w \preceq \gamma$ has length $k$, it follows that, for $i \neq j$,
$R(\epsilon_i)R(\delta_j) \preceq \gamma$ and that $\mu (\epsilon_i) \cdot \delta_j = 0$,
by Lemma \ref{lem:basic}.

If $w_i = R(\delta_1)
\cdots R(\delta_i)$, then $w_i(\delta_i) = -\epsilon_i$.  Since
$\epsilon_i, \delta_i \in \mM(w_i)$, it follows from Lemma \ref{lem:tau-rho}
that $\mu(\epsilon_i) \cdot
\delta_i = - \mu (w_i (\delta_i)) \cdot \delta_i = \mu (\delta_i)
\cdot \delta_i = 1$. For the last statement, note that any $\tau \in \Phi^+
(w)$ can be written as a nonnegative linear combination $\tau = a_1 \delta_1
+ \dots + a_k \delta_k$, so that $\mu(\epsilon_i) \cdot
\tau = a_i \ge 0$.
\end{proof}
\begin{lemma}
If $\tau \in \Phi^+ (w)$ satisfies $\tau < \epsilon_i$ then $\mu(\epsilon_i)
\cdot \tau = 0$ and $\epsilon_i \cdot \tau \ge 0$.
\label{lem:first1}
\end{lemma}
\begin{proof}
From Lemma \ref{lem:epsilon-delta} we have $\mu(\epsilon_i) \cdot \tau \ge 0$.
But $\tau < \epsilon_i$ implies $\mu(\epsilon_i) \cdot \tau \le 0$
by Lemma \ref{lem:BW2} (ii).
Hence $\mu(\epsilon_i) \cdot \tau = 0$ or, equivalently,
$R(\epsilon_i) R(\tau) \preceq \gamma$.
Finally $\epsilon_i \cdot \tau \ge 0$ by applying Lemma \ref{lem:BW3}
to this last relation and $\tau < \epsilon_i$.
\end{proof}
\begin{lemma}
Let $\tau \in \Phi^+ (w)$ and fix $1 \le i \le k$.  Then $\tau \in
\{\epsilon_1,\dots,\epsilon_i\}$ if and only if $R(\delta_i) R(\delta_{i-1})
\cdots R(\delta_1) \, \tau$ is a negative root.
In particular, $\tau \in \{\epsilon_1,\dots,\epsilon_k\}$ if and only if
$w^{-1} (\tau)$ is a negative root.
\label{lem:first0}
\end{lemma}
\begin{proof}
If $\tau = \epsilon_j$ with
$1 \le j \le i$ then $R(\delta_i) R(\delta_{i-1}) \cdots R(\delta_1) \,
\tau = - R(\delta_i) \cdots R(\delta_{j+1}) \, \delta_j$ is a
negative root since the $\delta_a \cdot \delta_b \le 0$.
For the converse let $\tau_0 = \tau$ and $\tau_j =
R(\delta_j) \tau_{j-1}$ for $1 \le j \le i$ and assume that $R(\delta_i)
R(\delta_{i-1}) \cdots R(\delta_1) \, \tau = \tau_i$ is a negative root.
Then there exists $1 \le j \le i$ such that $\tau_{j-1}$ is a positive
root but $\tau_j$ is negative. Since $\tau_j = R(\delta_j) \tau_{j-1}$
we must have $\tau_{j-1} = \delta_j$. On the other hand $\tau_{j-1} =
R(\delta_{j-1}) \cdots R(\delta_1) \, \tau$, so that $\tau = \epsilon_j$.
\end{proof}
\begin{lemma}
If $\tau \in \Phi^+ (w)$ satisfies $\tau < \epsilon_i$ and $\tau \not \in
\{\epsilon_1,\dots ,\epsilon_{i-1}\}$ then $\epsilon_i \cdot \tau = 0$.
\label{lem:first2}
\end{lemma}
\begin{proof}
From Lemma \ref{lem:first1} we have $\epsilon_i \cdot \tau \ge 0$ and
$\mu(\epsilon_i) \cdot \tau = 0$. Since $\tau \in \Phi^+ (w)$, the last
equation implies that $\tau$ lies in the positive hull of $\Pi (w) -
\{\delta_i\}$. Note that the root $\tau' = R(\delta_{i-1}) \cdots
R(\delta_1) \, \tau$, being positive by Lemma \ref{lem:first0}, lies
in the positive hull of $\Pi (w) - \{\delta_i\}$ as well. We conclude
that
\[ \epsilon_i \cdot \tau = (R(\delta_1) \cdots R(\delta_{i-1})
\, \delta_i) \cdot \tau = \delta_i \cdot \tau' \le 0 \]
and hence that $\epsilon_i \cdot \tau = 0$.
\end{proof}
\begin{corollary}
If $i < j$ and $\epsilon_i > \epsilon_j$ then $\epsilon_i \cdot \epsilon_j
= 0$.
\label{cor:first1}
\end{corollary}
\begin{proof}
Set $\tau = \epsilon_j$ in Lemma \ref{lem:first2}.
\end{proof}
%

%(The next statement will not be needed in the sequel but it is included here
%for reasons of completeness.
%
%\begin{corollary}
%If $i < j$ and $\epsilon_i > \epsilon_j$ then $\delta_i \cdot \delta_j = 0$.
%\label{cor:first2}
%\end{corollary}
%
%\begin{proof}
%Consider the root $R(\delta_{i+1}) \cdots R(\delta_{j-1}) \, \delta_j$. This
%can be expressed as a linear combination $a_{i+1} \delta_{i+1}+ \cdots +
%a_{j-1}\delta_{j-1} + a_j \delta_j$ with nonnegative coefficients, where
%$a_j = 1$. From Corollary \ref{cor:first2} we have
%
%\begin{eqnarray*}
%0 &=& \epsilon_i \cdot \epsilon_j\\
%&=& (R(\delta_1) \cdots R(\delta_{i-1}) \, \delta_i) \cdot
%(R(\delta_1) \cdots R(\delta_{j-1}) \, \delta_j)\\
%&=&  \delta_i \cdot (R(\delta_i) \cdots R(\delta_{j-1}) \, \delta_j)\\
%&=&  R(\delta_i) \, \delta_i \cdot (R(\delta_{i+1}) \cdots R(\delta_{j-1})
%\, \delta_j)\\
%&=& - \delta_i \cdot (R(\delta_{i+1}) \cdots R(\delta_{j-1}) \, \delta_j)\\
%&=& - \delta_i \cdot (a_{i+1}\delta_{i+1} + \cdots + a_{j-1} \delta_{j-1}
%+a_j \delta_j)\\
%&=& - a_{i+1}( \delta_i \cdot \delta_{i+1}) - \cdots - a_{j-1}( \delta_i
%\cdot \delta_{j-1}) - a_j ( \delta_i \cdot \delta_j). \\
%\end{eqnarray*}
%
%\vspace*{-0.1 in}
%\noindent
%Since the terms in the last expression are all nonnegative, for each $i <
%p \le j$ we must have either $a_p = 0$ or $\delta_i \cdot \delta_p = 0$.
%From $a_j = 1$ we deduce that $\delta_i \cdot \delta_j = 0$.
%\end{proof})
%
\begin{proposition}
The set $\{\epsilon_1,\dots,\epsilon_k\}$ is the vertex set of the first
facet of $X(w)$ in the lexicographic order.
\label{prop:first1}
\end{proposition}
\begin{proof}
It follows from (\ref{def:epsilon}) that $R (\epsilon_k) \cdots
R(\epsilon_1) = R(\delta_1) \cdots R(\delta_k) = w$. Moreover, by
Corollary \ref{cor:first1} we may relabel so that $\epsilon_1 < \dots < \epsilon_k$,
and hence the set
$\{\epsilon_1,\dots,\epsilon_k\}$ is the vertex set of some facet of
$X(w)$. Note that Lemma \ref{lem:first1} continues to hold for this reordered set. Given
$1 \le i \le k$ we have $\mu(\epsilon_i) \cdot \tau = 0$ for any $\tau
\in \Phi^+ (w)$ satisfying $\tau < \epsilon_i$ and, by our convention on
the $\epsilon$'s, we have $\mu(\epsilon_j) \cdot \tau = 0$ for $j \ge i$
as well. This forces such a root $\tau$ into a linear space of dimension
$i-1$, namely the intersection of $\mM(w)$ with the hyperplanes
$\mu(\epsilon_j)^\perp$ for $i \le j \le k$. Thus, for any facet $F$ of
$X(w)$, the $i$th vertex of $F$ cannot precede $\epsilon_i$ in the order
(\ref{total}).
\end{proof}
\begin{corollary}
Let $\tau \in \Phi^+ (w)$. Then $w^{-1} (\tau)$ is a negative root if and
only if $\tau$ is a vertex of the first facet of $X(w)$ in the
lexicographic order.
\label{cor:first3}
\end{corollary}
\begin{proof}
Combine Proposition \ref{prop:first1} with the last statement of Lemma
\ref{lem:first0}.
\end{proof}

\subsection{The lexicographically last facet}
\label{last}

We define the roots
\begin{equation}
\zeta_i = R(\delta_k) \cdots R(\delta_{i+1}) \, \delta_i
\label{def:zeta}
\end{equation}
for $1 \le i \le k$, so that $\zeta_k = \delta_k$. As with
$\{\epsilon_1,\dots,\epsilon_k\}$ in Section \ref{first} we see
that $\{\zeta_1,\dots,\zeta_k\}$ is a linearly independent subset
of $\Phi^+ (w)$. From (\ref{def:zeta}) we have
\begin{equation}
w = R(\delta_1) \cdots R(\delta_{i-1}) R(\delta_{i+1}) \cdots
R(\delta_k) R(\zeta_i) \label{eq:zeta2}
\end{equation}
which, combined with (\ref{eq:epsilon2}), gives $R(\epsilon_i) w =
w R(\zeta_i)$ or $R(\zeta_i)  =  w^{-1}R(\epsilon_i)w$ and hence
that $\zeta_i = \pm w^{-1} (\epsilon_i)$.
Since $w^{-1} (\epsilon_i)$ is a negative root by Lemma
\ref{lem:first0} we deduce that $\zeta_i = - w^{-1} (\epsilon_i)$.
\begin{lemma}
We have
\[ \mu(\delta_i) \cdot \zeta_j =
\begin{cases} 1 & \text{if \ $i=j$} \\
0 & \text{if \ $i\neq j$}. \end{cases} \]
In particular $\mu(\tau) \cdot \zeta_i \ge 0$ for all $i$ and for all
$\tau \in \Phi^+ (w)$. \label{lem:zeta-delta}
\end{lemma}
\begin{proof}
Using Lemma \ref{lem:tau-rho} (i) we get
\[ \mu(\delta_i) \cdot \zeta_j = -\mu(\delta_i) \cdot w^{-1}
(\epsilon_j) = - w(\mu(\delta_i)) \cdot \epsilon_j = \mu(\epsilon_j)
\cdot \delta_i \]
so the first statement follows from Lemma \ref{lem:epsilon-delta}. The
second statement follows as in the proof of Lemma \ref{lem:epsilon-delta}.
\end{proof}
The proofs of Lemmas~\ref{lem:last1}, \ref{lem:last0} and \ref{lem:last2} and of
Corollary~\ref{cor:last1} below
are completely analogous to
those of corresponding statements of Section \ref{first} and are omitted.
\begin{lemma}
If $\tau \in \Phi^+ (w)$ satisfies $\tau > \zeta_i$ then
$\mu(\tau) \cdot \zeta_i = 0$ and $\zeta_i \cdot \tau \ge 0$. \qed
\label{lem:last1}
\end{lemma}
\begin{lemma}
Let $\tau \in \Phi^+ (w)$. Given $1 \le i \le k$ we have $\tau \in
\{\zeta_i, \zeta_{i+1},\dots,\zeta_k\}$ if and only if
$R(\delta_i) R(\delta_{i+1}) \cdots R(\delta_k) \, \tau$ is a
negative root. In particular, $\tau \in \{\zeta_1,\dots,\zeta_k\}$
if and only if $w (\tau)$ is a negative root. \qed
\label{lem:last0}
\end{lemma}
\begin{lemma}
If $\tau \in \Phi^+ (w)$, $\tau > \zeta_i$ and $\tau \not \in
\{\zeta_{i+1},\dots,\zeta_k\}$ then $\zeta_i \cdot \tau  = 0$.
\qed \label{lem:last2}
\end{lemma}
\begin{corollary}
If $i < j$ and $\zeta_i > \zeta_j$ then $\zeta_i \cdot \zeta_j =
0$. \qed \label{cor:last1}
\end{corollary}
%
%\begin{corollary}
%If $i < j$ and $\zeta_i > \zeta_j$ then $\delta_i \cdot \delta_j =
%0$. \qed \label{cor:last2}
%\end{corollary}
%
\begin{proposition}
The set $\{\zeta_1,\dots,\zeta_k\}$ is the vertex set of the last
facet of $X(w)$ in the lexicographic order. \label{prop:last1}
\end{proposition}
\begin{proof}
As in the proof of Proposition \ref{prop:first1} we have $R
(\zeta_k) \cdots R(\zeta_1) = w$ and we may relabel so  that
$\zeta_1 < \dots < \zeta_k$ and hence the set
$\{\zeta_1,\dots,\zeta_k\}$ is
the vertex set of some facet of $X(w)$. Given $1 \le i \le k$,
Lemma \ref{lem:last1} implies that for $1 \le j \le i$ and any
$\tau \in \Phi^+ (w)$ with $\tau > \zeta_i$ we have $\mu(\tau)
\cdot \zeta_j = 0$. Since $\mu$ is an invertible linear
transformation, it follows as in the proof of Proposition
\ref{prop:first1} that any such root $\tau$ lies in a linear space
of dimension $k-i$ and hence that for any facet $F$ of $X(w)$, the
$i$th vertex of $F$ cannot succeed $\zeta_i$ in the order
(\ref{total}).
\end{proof}
\begin{corollary}
Let $\tau \in \Phi^+ (w)$. Then $w (\tau)$ is a negative root if and only
if $\tau$ is a vertex of the last facet of $X(w)$ in the lexicographic
order.
\label{cor:last3}
\end{corollary}
\begin{proof}
Combine Proposition \ref{prop:last1} with the last statement of Lemma
\ref{lem:last0}.
\end{proof}

\begin{remark}
{\rm In the special case $w = \gamma$ the vertices
of the last facet of $X(\gamma)$ are the last $n$ positive roots
by Remark~\ref{rem:lastn}.
}
\label{rem: Omega}
\end{remark}

\begin{remark}
{\rm As in the proofs of Propositions \ref{prop:first1} and
\ref{prop:last1}, in the following sections we will denote by
$\{\epsilon_1, \epsilon_2,\dots,\epsilon_k\}$ and $\{\zeta_1,
\zeta_2,\dots,\zeta_k\}$ the ordered vertex sets of the first and
last facet of $X(w)$, respectively, in the lexicographic order.
Clearly Lemmas \ref{lem:first1}, \ref{lem:first2}, \ref{lem:last1}
and \ref{lem:last2}, as well as the last statement of Lemmas
\ref{lem:epsilon-delta} and \ref{lem:zeta-delta}, continue to hold
under this convention.  } \label{rem:ezeta}
\end{remark}

\begin{remark}
{\rm It follows from Proposition~\ref{prop:last1} that, for $1 < i \le k$,
the set $\{\zeta_i,\dots,\zeta_k\}$ is the vertex set of the last facet of
$X(wR(\zeta_1) \cdots R(\zeta_{i-1}))$.
Similarly, it follows from Proposition~\ref{prop:first1} that, for $1 \le i
< k$, the set $\{\epsilon_1,\dots,\epsilon_i\}$ is the vertex set of the first
facet of $X(R(\epsilon_{i+1}) \cdots R(\epsilon_k)w )$.
} \label{rem:facetfaces}
\end{remark}

\section{Vertex type}
\label{sec:sub}

In this section $F$ denotes a face of $\Delta (\gamma)$ with
ordered vertex set $\{\tau_1, \tau_2,\dots,\tau_k\}$, $w =
R(\tau_k) R(\tau_{k-1}) \cdots R(\tau_1) \preceq \gamma$ and
$\{\epsilon_1, \epsilon_2,\dots,\epsilon_k\}$ and $\{\zeta_1,
\zeta_2,\dots,\zeta_k\}$ are the ordered vertex sets of the first
and last facet of $X(w)$, respectively, in the lexicographic
order.
\begin{definition}
%%Let $V = \{\tau_1, \tau_2,\dots,\tau_k\}$ be an ordered subset of the
%%vertex set of $\Delta (\gamma)$ such that
%
%%\[ R(\tau_k) R(\tau_{k-1}) \cdots R(\tau_1) \preceq \gamma. \]
%
For $1 \le i \le k$ we set
\[ \begin{tabular}{rcl}
$u_i (F)$ & $\! \! =$ & $\! \! R(\tau_k) R(\tau_{k-1}) \cdots R(\tau_i)$ \\
$v_i (F)$ & $\! \! =$ & $\! \! R(\tau_1) R(\tau_2) \cdots R(\tau_i)$
\end{tabular} \]
and say that $\tau_i$ is a left vertex in $F$ if $u_i (F) \tau_i$ is a
positive root. Otherwise we say that $\tau_i$ is a right vertex in $F$.
\label{def:type}
\end{definition}
%
%%Note that the order on $V$ need not be consistent with the total order on
%%the roots. However we will be mostly concerned with the special case in
%%which $V$ is consistently ordered, so that $V$ is the vertex set of some
%%face $F$ of $GA(W)$. In this case we will write $u_i (F)$ and $v_i (F)$
%%instead of $u_i (V)$ and $v_i (V)$, respectively, and will refer to left
%%and right vertices in $V$ as left and right vertices in $F$.

%%For the remainder of this section $F$ denotes a face of $GA(W)$ with
%%ordered vertex set $\{\tau_1, \tau_2,\dots,\tau_k\}$ and $w = R(\tau_k)
%%R(\tau_{k-1}) \cdots R(\tau_1) \preceq \gamma$.
%
\begin{lemma}
Any vertex of $F$ which is a negative root is a left vertex in $F$.
\label{lem:easyleft}
\end{lemma}
\begin{proof}
Suppose that $\tau_i \in -\Pi$ and recall that the first $n-s$
vertices of $\Delta (\gamma)$ in the total order (\ref{total}) are
in $-\Pi_2$ while the last $s$ are in $-\Pi_1$. If $\tau_i \in
-\Pi_2$ then $\tau_1,\dots,\tau_{i-1} \in -\Pi_2$ as well. In
particular $\tau_i$ is orthogonal to each of
$\tau_1,\dots,\tau_{i-1}$ and hence $u_i (F) \, \tau_i = w
(\tau_i)$. Since $-\tau_i \in \Pi_2 \subseteq \Omega$, Corollary
\ref{lem:nsteps} (iii) implies that $w (\tau_i)$ is a positive
root, so that $\tau_i$ is a left vertex in $F$. Similarly, if
$\tau_i \in -\Pi_1$ then $\tau_{i+1},\dots,\tau_{n} \in -\Pi_1$.
In particular $\tau_i$ is orthogonal to each of
$\tau_{i+1},\dots,\tau_k$ and hence $u_i (F) \, \tau_i = R(\tau_i)
\tau_i = -\tau_i$ is a positive root, so that $\tau_i$ is a left
vertex in $F$.
\end{proof}
\begin{lemma}
The vertex $\tau_i$ is a right vertex in $F$ if and only if $\tau_i$ is
a vertex of the last facet of $X(u_i(F))$ in the lexicographic order.
Moreover, if $F$ is positive then $\tau_i$ is the first vertex of this
facet.
\label{lem:type1}
\end{lemma}
\begin{proof}
Since any right vertex in $F$ must be a positive root by Lemma
\ref{lem:easyleft}, the first statement is a direct consequence of the
definition and Corollary \ref{cor:last3}. The second statement is obvious
since $\{\tau_i,\dots,\tau_k\}$ is the ordered vertex set of a facet of
$X(u_i(F))$.
\end{proof}
Let $F$ be positive, so that $F$ is a facet of $X(w)$. We will
describe the wall of $F$ opposite (not containing) the vertex
$\tau_i$. To simplify notation, define the roots $\eta_i = u_i (F) \tau_i$ and
$\theta_i = v_i (F) \tau_i$.   Then
\[ R(\eta_i) R(\tau_k) \cdots R(\tau_{i+1}) R(\tau_{i-1}) \cdots
R(\tau_1) = w \preceq \gamma \]
and
\[ R(\tau_k) \cdots R(\tau_{i+1}) R(\tau_{i-1}) \cdots R(\tau_1) R(\theta_i)
= w \preceq \gamma. \]
Note that $\eta_i, \theta_i \in \Phi (w)$ for all $i$, although they are
not necessarily positive roots.  As
in the proof of Lemma \ref{lem:epsilon-delta} we find that
\begin{equation}
\mu(\eta_i) \cdot \tau_j =
\left \{
\begin{array}{rr}
-1, & \text{if \ $i=j$} \\
0, & \text{if \ $i\neq j$,}
\end{array}\right.
\label{eq:betas}
\end{equation}
\begin{equation}
\mu(\tau_j) \cdot \theta_i =
\left \{
\begin{array}{rr}
-1, & \text{if \ $i=j$} \\
0, & \text{if \ $i\neq j$.}
\end{array}\right.
\label{eq:beta's}
\end{equation}
As a result, the wall of $F$ opposite $\tau_i$ is the intersection
of $S^{n-1} \cap \mM(w)$ with $\mu (\eta_i)^{\perp}$, the linear
hyperplane orthogonal to
$\mu (\eta_i)$. It follows by linearity from (\ref{eq:betas}) and
(\ref{eq:beta's}) that $\mu(\eta_i) \cdot x = \mu(x) \cdot
\theta_i$ for all $x \in \mM(w)$. Note that, since $X(w)$ is
homeomorphic to a ball, every codimension one face is contained in
either exactly one or exactly two facets of $X(w)$ and the first
case occurs if and only if this face is contained in the boundary
of $X(w)$.
\begin{proposition}
Suppose that $F$ is positive.
\begin{enumerate}
\itemsep=0pt
\item[(i)] The root $\tau_i$ is a left vertex in $F$ if and only if there
exists a vertex $\tau > \tau_i$ of $X(w)$ such that $(V \sm \, \{\tau_i\})
\cup \{\tau\}$ is the vertex set of a facet of $X(w)$. Moreover, such a
vertex $\tau$ is unique.

\item[(ii)] The root $v_i (F) \tau_i$ is positive if and only if there
exists a vertex $\tau < \tau_i$ of $X(w)$ such that $(V \sm \, \{\tau_i\})
\cup \{\tau\}$ is the vertex set of a facet of $X(w)$. Moreover, such a
vertex $\tau$ is unique.
\end{enumerate}
\label{prop:type2}
\end{proposition}
\begin{proof}
The uniqueness of $\tau$ in both parts follows from the discussion preceding the
statement of the proposition.

\medskip
\noindent (i) Suppose $\tau_i$ is a left vertex in $F$, so that $\eta_i = u_i (F)
\tau_i$ is a positive root.  Since $\mu (\eta_i) \cdot \tau_i < 0$ by
(\ref{eq:betas}) and $\mu (\eta_i) \cdot \eta_i = 1$ we see that
$\mu (\eta_i)^{\perp}$ separates the positive roots $\tau_i$ and $\eta_i$. It follows
that the wall of $F$ opposite $\tau_i$ does not lie in the boundary of $X(w)$
and, as a result, there exists a facet of $X(w)$ other than $F$ with vertex set
$(V \sm \, \{\tau_i\}) \cup \{\tau\}$.  Since
$\mu (\eta_i)^{\perp}$ separates $\tau$ from $\tau_i$, we must have
$\mu (\eta_i) \cdot \tau > 0$. In
view of Lemma \ref{lem:BW2}, this inequality combined with $\mu (\eta_i)
\cdot \tau_i < 0$ gives $\tau > \tau_i$.  For the converse suppose that
$(V \sm \, \{\tau_i\}) \cup \{\tau\}$ is the vertex set of a facet of
$X(w)$ for some $\tau > \tau_i$. Then $\mu (\eta_i) \cdot \tau_i < 0$
and $\mu (\eta_i) \cdot \tau > 0$, hence Lemma \ref{lem:BW2} implies
that $\eta_i$ is a positive root, meaning that $\tau_i$ is a left vertex
in $F$.

\medskip
\noindent (ii) Suppose that $\theta_i = v_i (F) \tau_i$ is a positive root.
Since $\mu (\eta_i) \cdot \tau_i < 0$ by
(\ref{eq:betas}) and $\mu (\eta_i) \cdot \theta_i = \mu (\theta_i) \cdot \theta_i = 1$
we see that $\mu (\eta_i)^{\perp}$ separates the positive roots $\tau_i$ and $\theta_i$.
As in part (i), it follows that the wall of $F$ opposite $\tau_i$
does not lie in the boundary of $X(w)$ and, as a result, there exists
a facet of $X(w)$ other than $F$ with vertex set $(V \sm \, \{\tau_i\})
\cup \{\tau\}$.  Since $\mu (\eta_i)^{\perp}$ separates $\tau$
from $\tau_i$, we must have
$\mu (\eta_i) \cdot \tau > 0$.
Thus $\mu (\tau) \cdot \theta_i  = \mu (\eta_i) \cdot \tau > 0$.
In view of Lemma \ref{lem:BW2},
this inequality combined with
$\mu (\tau_i) \cdot \theta_i = \mu (\eta_i) \cdot \tau_i < 0$
gives $\tau < \tau_i$. The proof of the converse proceeds as in
part (i).
\end{proof}
\begin{corollary}
If $F$ is a positive face and $\tau_i$ is a left vertex in $F$ then
$v_i (F) \tau_i$ is a negative root.
\label{lem:rho-left}
\end{corollary}
\begin{proof}
By part (i) of Proposition~\ref{prop:type2}, there
exists a vertex $\tau > \tau_i$ of $X(w)$ such that $(V \sm \, \{\tau_i\})
\cup \{\tau\}$ is the vertex set of a facet of $X(w)$.  Since $X(w)$ is a manifold
there cannot be a vertex $\tau' < \tau_i$ of $X(w)$ such that $(V \sm \, \{\tau_i\})
\cup \{\tau'\}$ is the vertex set of a facet of $X(w)$.  Thus the root $v_i (F) \tau_i$
cannot be positive by part (ii) of Proposition~\ref{prop:type2}.
%To give an alternative proof, let us abbreviate
%$u_i (F)$ and $v_i (F)$ as $u_i$ and $v_i$, respectively, so that
%$\theta_i = v_i (\tau_i)$. We have $\gamma (\theta_i) = \gamma v_i
%(\tau_i) = u_i R(\tau_i) (\tau_i) = -u_i (\tau_i)$. By assumption
%$u_i (\tau_i)$ is a positive root, so $\gamma(\theta_i)$ is a
%negative root. This implies that either $\theta_i$ is a negative
%root or that $\theta_i \in \Omega$. However if the later
%were true then $v_i^{-1} (\theta_i) = \tau_i$ is a positive root with
%$v_i^{-1} \preceq \gamma$, contradicting part (iii) of Corollary
%\ref{lem:nsteps}. Hence $\theta_i$ must be negative.
\end{proof}

\begin{corollary}
If $\tau \in \Phi^+ (w)$ satisfies $\tau \le \zeta_1$ then there
exists a facet of $X(w)$ having $\tau$ as its smallest vertex. In
particular, $\tau$ precedes all vertices of the last facet of $X(w
R(\tau))$ in the lexicographic order. \label{cor:tau}
\end{corollary}
\begin{proof}
The statement is obvious in case $\tau = \zeta_1$, so suppose
$\tau < \zeta_1$. Let $F$ be any facet of $X(w)$ having $\tau$ as
a vertex and let $\{\tau_1, \tau_2,\dots,\tau_k\}$ be the ordered
vertex set of $F$. If $\tau = \tau_1$ there is nothing to prove.
Otherwise, since $\tau_1 < \zeta_1$, the vertex $\tau_1$ is a left
vertex in $F$ by Corollary \ref{cor:last3}. Proposition
\ref{prop:type2} (i) implies that there exists a vertex $\tau'$ of
$X(w)$ with $\tau_1 < \tau'$ such that $(V \sm \, \{\tau_1\}) \cup
\{\tau'\}$ is the vertex set of a facet $F'$ of $X(w)$. Clearly
$\tau$ is a vertex of $F'$ and the smallest vertex of $F'$
succeeds that of $F$ in the order (\ref{total}). Therefore,
applying the same argument to $F'$ repeatedly, if necessary, we
can find a facet of $X(w)$ having $\tau$ as its smallest vertex.
\end{proof}

The following technical fact will be used in Section \ref{bij}.
\begin{proposition}
If $\tau \in \Phi^+ (w)$ satisfies
\[\zeta_1 < \cdots < \zeta_{i-1} < \tau < \zeta_i < \cdots < \zeta_k \]
for some $1 \le i \le k$ then
\begin{enumerate}
\itemsep=0pt
\item[(i)] $\tau$ is orthogonal to $\zeta_1,\dots,\zeta_{i-1}$,
\item[(ii)] $R(\tau) \preceq w R(\zeta_1) \cdots R(\zeta_{i-1})$ and
\item[(iii)] the last facet of $X(w R(\tau))$ in the lexicographic
order has ordered vertex set
\[ \{\zeta_1,\dots,\zeta_{i-1}, \zeta'_{i+1},\dots,\zeta'_k\},\]
where $\{\zeta'_{i+1},\dots,\zeta'_k\}$ is the ordered vertex set
of the lexicographically last facet of $X(w R(\zeta_1) \cdots
R(\zeta_{i-1}) R(\tau))$. Moreover $\tau < \zeta'_{i+1}$.
\end{enumerate}
\label{prop:last2}
\end{proposition}
\begin{proof}
Part (i) follows from Lemma \ref{lem:last2}. In view of Lemma
\ref{lem:BW} and the fact that $R(\zeta_{i-1}) \cdots R(\zeta_1)
\preceq \gamma$ has length $i-1$, we can deduce from Lemma
\ref{lem:last1} that $R(\tau) R(\zeta_{i-1}) \cdots R(\zeta_1)
\preceq \gamma$ has length $i$. Since $R(\tau) \preceq w$ and
$R(\zeta_{i-1}) \cdots R(\zeta_1) \preceq w$ by the assumptions,
it follows from Lemma \ref{lem:basic1} (iv) that
\[ R(\tau) R(\zeta_{i-1}) \cdots R(\zeta_1) \preceq w, \]
which proves (ii). Let $w' = w R(\zeta_1) \cdots R(\zeta_{i-1}) =
R(\zeta_k) \cdots R(\zeta_i)$ and observe that
$\{\zeta_i,\dots,\zeta_k\}$ is the ordered vertex set of the last
facet of $X(w')$ in the lexicographic order by Remark~\ref{rem:facetfaces}.
Since $R(\tau)
\preceq w'$ and $\tau < \zeta_i$, Corollary \ref{cor:tau} implies
that $\tau < \zeta'_{i+1}$ and, in particular, $\zeta_{i-1} <
\zeta'_{i+1}$. Since $R(\tau)$ commutes with $R(\zeta_j)$ for $1
\le j \le i-1$ by (i), we have
\[ R(\zeta'_k) \cdots R(\zeta'_{i+1}) = w' R(\tau) = w R(\tau)
R(\zeta_1) \cdots R(\zeta_{i-1}) \]
and hence the set $\{\zeta_1,\dots,\zeta_{i-1},
\zeta'_{i+1},\dots,\zeta'_k\}$ is the ordered vertex set of a
facet of $X(w R(\tau))$. To complete the proof of (iii) it
suffices to show that $\zeta_j$ is the $j$th vertex of the last
facet of $X(w R(\tau))$ in the lexicographic order for $1 \le j <
i$. This follows from the claim that any positive root $\tau' >
\zeta_j$ in $\mM (w R(\tau))$ lies in the $(k-j-1)$-dimensional
space $\mM (w R(\zeta_1) \cdots R(\zeta_j) R(\tau))$. Indeed, we
have $R(\tau') \preceq w R(\tau)$, $R(\zeta_j) \cdots R(\zeta_1)
\preceq w R(\tau)$ and, by Lemmas \ref{lem:BW} and
\ref{lem:last1}, $R(\tau') R(\zeta_j) \cdots R(\zeta_1) \preceq
\gamma$. From Lemma \ref{lem:basic1} (iv) we deduce that $R(\tau')
R(\zeta_j) \cdots R(\zeta_1) \preceq w R(\tau)$.
\end{proof}

\section{The map $\phi$}
\label{types}

\begin{definition}
For a facet $F$ of $\Delta (\gamma)$ with ordered vertex set $\{\sigma_1,
\sigma_2,\dots,\sigma_n\}$ we define
\[ \phi(F) = c(F, \sigma_n) c(F, \sigma_{n-1}) \cdots c(F, \sigma_1) \]
where
\[ c(F, \sigma_i) = \begin{cases}
R(\sigma_i), & \text{if $\sigma_i$ is a right vertex in $F$} \\
I, & \text{otherwise.}
\end{cases} \]
\label{def:phi}
\end{definition}

\medskip
In Figure \ref{f:cyclohedron}, the right vertices
of each facet are indicated with a small circle and the value of
$\phi$ on each facet can be deduced.  The vertex types can be verified using
Definition~\ref{def:type}, but the calculations are simplified greatly
by using Lemma~\ref{lem:rho-left}, Theorem~\ref{thm:lexlast} and
Lemma~\ref{lem:neg-simple}.
\medskip

Clearly $\phi(F) \in \nc (\gamma)$ and thus $\phi$ is a map from
the set of facets of $\Delta (\gamma)$ to $\nc (\gamma)$ such
that, for any facet $F$ of $\Delta (\gamma)$, the rank of
$\phi(F)$ in $\nc (\gamma)$ is equal to the number of right
vertices of $F$.
 It will be shown in Section \ref{bij} that $\phi$ is
a bijection. Part of the injectivity of $\phi$ will be proved in
this section. The next theorem gives a fundamental property of
$\phi$.
\begin{theorem}
Let $F$ be a positive facet of $\Delta (\gamma)$ and let $w = \phi(F)$.
\begin{enumerate}
\itemsep=0pt
\item[(i)] The set of right vertices in $F$ is equal to the vertex set
of the last facet of $X(w)$ in the lexicographic order.
\item[(ii)] The set of left vertices in $F$ is equal to the vertex set
of the first facet of $X(w^{-1} \gamma)$ in the lexicographic order.
\end{enumerate}
\label{thm:lexlast}
\end{theorem}

For the proof of the theorem we need the following lemma.
\begin{lemma}
Suppose that $F$ is a positive face of $\Delta (\gamma)$ with ordered vertex
set $\{\tau_1, \tau_2,\dots,\tau_k\}$. If $\tau_i$ is a right vertex and
$\tau_j$ is a left vertex in $F$ for some $i<j$ then $\tau_i \cdot \tau_j
= 0$.
\label{lem:commute}
\end{lemma}
\begin{proof}
Proceeding by induction on $j-i$, we may assume that the result
holds (i) for the face with ordered vertex set
$\{\tau_i,\dots,\tau_{j-1}\}$ and (ii) for the face with ordered
vertex set $\{\tau_{i+1},\dots,\tau_j\}$. Furthermore, we may also
assume that all roots in $\{\tau_{i+1},\dots,\tau_{j-1}\}$ are
right vertices in $F$. Indeed, if not then we replace $F$ by the
face $F'$ obtained from $F$ by removing all roots in
$\{\tau_{i+1},\dots,\tau_{j-1}\}$ which are left vertices in $F$
and observe that, in view of (i), $\tau_i$ and $\tau_j$ are still
right and left vertices in $F'$, respectively. Under these
assumptions, let $w = u_i (F)$ and $\{\zeta_1,
\zeta_2,\dots,\zeta_r\}$ be the ordered last facet of $X(w)$ in
the lexicographic order, where $r = k-i+1$. We have $\tau_i =
\zeta_1$ by Lemma \ref{lem:type1} and
\[ R(\tau_j) \preceq u_{i+1} (F) = w R(\tau_i) = R(\zeta_r)
R(\zeta_{r-1}) \cdots R(\zeta_2). \]
By Remark~\ref{rem:facetfaces}, the last facet of $X(u_{i+1} (F))$
in the lexicographic order has vertex set
$\{\zeta_2,\dots,\zeta_r\}$. By our assumptions $\tau_j$ is
orthogonal to all roots in $\{\tau_{i+1},\dots,\tau_{j-1}\}$ and
hence $u_{i+1} (F) \tau_j = u_j (F) \tau_j$. Since $\tau_j$ is a
left vertex in $F$ we conclude that $u_{i+1} (F) \tau_j$ is a
positive root. It follows from Lemma \ref{lem:last0} that $\tau_j
\not \in \{\zeta_2,\dots,\zeta_r\}$. Since $\tau_j \in \mM(w)$ and
$\tau_j > \zeta_1 = \tau_i$, Lemma \ref{lem:last2} applies to give
$\tau_i \cdot \tau_j = 0$.
\end{proof}

\bigskip
\noindent
\emph{Proof of Theorem \ref{thm:lexlast}.}
Let $\{\sigma_1, \sigma_2,\dots,\sigma_n\}$ be the ordered vertex set of
$F$ and $\{\sigma_{i_1}, \sigma_{i_2},\dots,\sigma_{i_k}\}$ and
$\{\sigma_{j_1}, \sigma_{j_2},\dots,\sigma_{j_\ell}\}$ be the ordered
sets of right and left vertices in $F$, respectively, so that $w =
R(\sigma_{i_k}) \cdots R(\sigma_{i_1})$. To prove (i) it suffices to show
that $\sigma_{i_p}$ is the first vertex of the last facet of
$X(R(\sigma_{i_k}) \cdots R(\sigma_{i_p}))$ in the lexicographic order
for $1 \le p \le k$. In view of Corollary \ref{cor:last3}, this is equivalent
to the statement that the root $R(\sigma_{i_k}) \cdots R(\sigma_{i_p})
\sigma_{i_p}$ is negative. The last statement holds since $R(\sigma_{i_k})
\cdots R(\sigma_{i_p}) \sigma_{i_p} = u_{i_p} (F) \sigma_{i_p}$ by Lemma
\ref{lem:commute} and $\sigma_{i_p}$ is a right vertex in $F$.

Similarly, by Lemma \ref{lem:commute} we have $w^{-1} \gamma =
R(\sigma_{j_\ell}) \cdots R(\sigma_{j_1})$. To prove (ii) it suffices
to show that $\sigma_{j_q}$ is the last vertex of the lexicographically
first facet of $X(R(\sigma_{j_q}) \cdots R(\sigma_{j_1}))$ for $1 \le q \le
\ell$. By Corollary \ref{cor:first3} this is equivalent to the statement
that $R(\sigma_{j_1}) \cdots R(\sigma_{j_q}) \sigma_{j_q}$ is a negative
root. Since $R(\sigma_{j_1}) \cdots R(\sigma_{j_q}) \sigma_{j_q} = v_{j_q}
(F) \sigma_{j_q}$ by Lemma \ref{lem:commute}, this follows from Corollary
\ref{lem:rho-left}.
\qed

\medskip
\begin{corollary}
The restriction of the map $\phi$ to the set of positive facets of
$\Delta (\gamma)$ is injective.
\label{cor:inj}
\end{corollary}
\begin{proof}
If $F$ is a positive facet of $\Delta (\gamma)$ with $w = \phi(F)$, then
parts (i) and (ii) of Theorem~\ref{thm:lexlast} imply that the vertices of $F$ are
determined by $w$ and hence $F$ is the unique positive facet whose image under $\phi$
is $w$.
\end{proof}
The image of the facets of $X(\gamma)$ can be characterized in terms of non-peripheral
elements. Recall from Section~\ref{peripheral} that $w \in L(\gamma)$ is
peripheral if and only if  $R(\rho) \preceq w^{-1}\gamma$ for some $\rho
\in \Omega$.
\begin{proposition}
If $F$ is a facet of $\Delta (\gamma)$ then the element $\phi(F)$ is non-peripheral
if and only if $F$ is positive.
\label{prop:per}
\end{proposition}
\begin{proof}
Let $\{\sigma_1, \sigma_2,\dots,\sigma_n\}$ be the ordered vertex set of
$F$. Assume that $\sigma_i$ is a negative simple root for some $1 \le
i \le n$, say $\sigma_i = -\alpha_j \in -\Pi_1$ with $1 \le j \le s$ (the
case $\sigma_i \in -\Pi_2$ is similar). As in proof of Lemma \ref{lem:easyleft},
$R(\sigma_i)$ commutes with each of $R(\sigma_{i+1}),\dots,R(\sigma_n)$.
By the same lemma $\sigma_i$ is a left vertex in $F$ and hence $$\phi(F)
\preceq R(\sigma_n) \cdots R(\sigma_{i+1}) R(\sigma_{i-1}) \cdots
R(\sigma_1) = R(\sigma_i) \gamma = R(\alpha_j) \gamma,$$ which implies that
$\phi(F)$ is peripheral.

The converse is proved by contradiction. Let $F$ be positive and assume
that $\phi(F)$ is peripheral.  Set $w =
\phi(F)$, $v = w^{-1} \gamma$ and $r = \ell(v)$. Let
$\{\delta'_1,\dots,\delta'_r\}$ be the ordered simple system for
$\Phi^+(v)$ and let $\{\epsilon'_1,\dots,\epsilon'_r\}$ be the set
of vertices of the first facet of $X(v)$, defined as in
equation~(\ref{def:epsilon}). Recall from \cite[Theorem 5.1]{BW}
that $\delta'_r$ is the largest root in $\Phi^+ (v)$ with respect
to the order (\ref{total}). Since $w$ is peripheral, the set
$\Omega$ intersects $\Phi^+(v)$ by
Proposition~\ref{prop:per-char} and hence $\delta'_r \in
\Omega$. Since $\epsilon'_r$ is a vertex of the first facet of
$X(v)$ in the lexicographic order, $\epsilon'_r$ is a left vertex
in $F$ by Theorem \ref{thm:lexlast} (ii), say $\epsilon'_r =
\sigma_i$. Let $u_i = u_i (F) = R(\sigma_n) \cdots R(\sigma_i)$
and $v_i = v_i (F) = R(\sigma_1) \cdots R(\sigma_i)$.
By Corollary~\ref{cor:first1} and Proposition~\ref{prop:first1}, the
root $\sigma_i$ is orthogonal to all later left vertices in $F$ (if any).
It follows from this fact and Lemma \ref{lem:commute} that
$v_i(\sigma_i) = v^{-1}(\sigma_i)$.
However $v^{-1}(\sigma_i) =
v^{-1}(\epsilon'_r) = -\delta'_r \in -\Omega$ and hence $\gamma
(v_i (\sigma_i)) \in -\gamma (\Omega)$ is a positive root by
Remark \ref{rem:lastn}. On the other hand,
$\gamma (v_i (\sigma_i)) = -u_i(\sigma_i)$ since $\gamma =
u_iR(\sigma_i)v_i^{-1}$.  But $ -u_i(\sigma_i)$ (and hence $\gamma (v_i (\sigma_i))$)
is a negative root since $\sigma_i$ is a left vertex in $F$.  This gives the
required contradiction.
\end{proof}

\section{Bijectivity of $\phi$}
\label{bij}

The first step in establishing bijectivity of $\phi$ is to show that
$\phi$ maps the set of facets of $X (\gamma)$ surjectively onto the set
of non-peripheral elements of $\nc (\gamma)$.

\begin{lemma}
If $w \in \nc (\gamma)$ is non-peripheral then there exists a facet $F$ of
$X(\gamma)$ such that the left vertices in $F$ are precisely the vertices
of the first facet of $X(w^{-1} \gamma)$ in the lexicographic order. In
particular, $\phi(F) = w$.
\label{lem:surj}
\end{lemma}
\begin{proof}
Let $\{\tau_1,\dots,\tau_k\}$ be the ordered vertex set of the first facet
of $X(w^{-1} \gamma)$ in the lexicographic order. For $0 \le i \le k$
we define
\[ v_i = R(\tau_i) \cdots R(\tau_1) \,\,\, \mbox{ and } \,\,\,
w_i = \gamma v_i^{-1}, \]
so that $v_0 = I$, $w_0 = \gamma$, $v_k = w^{-1}
\gamma$ and $w_k = w$. Note that $v_i = R(\tau_i) v_{i-1}$ and $w_{i-1}
= w_i R(\tau_i)$ for $1 \le i \le k$. We first claim that $w_{i-1}
(\tau_i)$ is a positive root for $1 \le i \le k$.  Since $\tau_i$
is a vertex of the first facet of $X(v_i)$ in the lexicographic order, the
root $v_i^{-1} (\tau_i)$ is negative by Corollary \ref{cor:first3}.
Therefore $-v_i^{-1} (\tau_i)$ is a positive root, clearly in $\mM(v_i)$.
Since $w$, and hence $w_i$, is non-peripheral we know from Proposition
\ref{prop:per-char} that this root cannot be in the set $\Omega$ of the
last $n$ positive roots. Hence $\gamma(-v_i^{-1}(\tau_i))$ must be a
positive root. This proves the claim since $w_{i-1}
(\tau_i) = - w_i (\tau_i) = - \gamma (v_i^{-1} (\tau_i))$.

We will show by induction that for each $0 \le i \le k$ there
exists a facet $F_i$ of $X(\gamma)$ such that the ordered set of
left vertices in $F_i$ is equal to $\{\tau_1,\dots,\tau_i\}$. For
$i = 0$ observe that the last facet $F_0$ of $X(\gamma)$ in the
lexicographic order, having $\Omega$ as its vertex set, has the
desired property by Lemma \ref{lem:type1}.
\vskip .2cm
For the
inductive step let $1 \le i \le k$ and assume that $X(\gamma)$ has
a facet $F_{i-1}$ whose ordered set of left vertices is
$\{\tau_1,\dots,\tau_{i-1}\}$. Let $\{\zeta_i,\dots,\zeta_n\}$ be
the ordered set of right vertices in $F_{i-1}$. Since $\gamma =
w_{i-1} v_{i-1}$, we know from Lemma \ref{lem:commute} and Theorem
\ref{thm:lexlast} (ii) that $\phi(F_{i-1}) = w_{i-1}$ and that
$\{\zeta_i,\dots,\zeta_n\}$ is the vertex set of the last facet of
$X(w_{i-1})$ in the lexicographic order. From $w_{i-1} = w_i
R(\tau_i)$ we get $R(\tau_i) \preceq w_{i-1}$.   Furthermore,
$\tau_i \notin \{\zeta_i,\dots,\zeta_n\}$ by Corollary \ref{cor:last3},
 since $w_{i-1}(\tau_i)$ is a positive root by the earlier claim. Thus
$\zeta_i < \dots < \zeta_{j-1} < \tau_i < \zeta_j < \dots <
\zeta_n$ for some $j$. Therefore Proposition \ref{prop:last2}
applies to $w_{i-1}$ and $\tau_i$ to give
\begin{enumerate}
\itemsep=0pt \item[(i)] $\tau_i$ is orthogonal to each of
$\zeta_i, \dots, \zeta_{j-1}$,
\item[(ii)] $R(\tau_i) \preceq
w_{i-1} R(\zeta_i)\dots R(\zeta_{j-1})$,
\item[(iii)] the ordered vertex set
of $G_i$, the last facet of $X(w_{i-1} R(\tau_i))$ in the
lexicographic order, is $\{\zeta_i,\dots,\zeta_{j-1},
\zeta'_{j+1},\dots,\zeta'_n\}$ and
\item[(iv)] $\tau_i < \zeta'_{j+1}$.
\end{enumerate}
Here $\{\zeta'_{j+1},\dots,\zeta'_n\}$ is the ordered vertex set of the
last facet of $X(w')$ in the lexicographic order, where
$w' = w_{i-1}R(\zeta_i) \dots R(\zeta_{j-1})R(\tau_i)$.
\vskip .2cm
The union of $\{\tau_1,\dots,\tau_{i-1}, \tau_i\}$ with
the vertex set of $G_i$ is the vertex set of a facet $F_i$ of
$X(\gamma)$, since the roots in
\[V = \{\tau_1,\dots,\tau_{i-1}\}\cup \{\zeta_i,\dots,\zeta_{j-1}\}\]
can be ordered as they were in $F_{i-1}$, while $\tau_i$ can be
positioned between $\zeta_{j-1}$ and $\zeta'_{j+1}$ by (iv). Since
the roots in $V$ are the first $j-1$ roots in both $F_{i-1}$ and
$F_i$ and are ordered in precisely the same way, it follows that
$\tau_r$ remains a left vertex in $F_i$ for all $1 \le r \le i-1$.
By (i) and (iii) $\zeta_t$ remains a right vertex in $F_i$ for all
$i \le t \le j-1$ and each $\zeta'_t$ for $j < t \le n$ is a right
vertex in $F_i$. Finally, using (i), we have that $w' R(\tau_i)
(\tau_i) = w_{i-1} R(\zeta_i) \cdots R(\zeta_{j-1}) (\tau_i) =
w_{i-1} (\tau_i)$ is a positive root by our claim, so that
$\tau_i$ is a negative vertex in $F_i$, as desired. Thus
$\phi(F_i) = w_i$ and the induction is complete.
\end{proof}

Let us call an element of $\nc (\gamma)$ of the from
\begin{equation}
\gamma' = R(\alpha_{i_1}) \cdots R(\alpha_{i_t}),
\label{gamma'}
\end{equation}
with $1 \le i_1 < \cdots < i_t \le n$, a \emph{standard parabolic Coxeter
element} of $W$ (with respect to the ordered simple system $\Pi$).
Note that the intersection of the moved spaces of two standard parabolic
Coxeter elements is again the moved space of some standard parabolic Coxeter
element. It follows from Lemma \ref{lem:basic1} (i) that, given $w \preceq
\gamma$, there is a minimum, with respect to the partial order $\preceq$,
standard parabolic Coxeter element $\gamma'$ satisfying $w \preceq \gamma'$.
Clearly $w$ is non-peripheral with respect to $\gamma'$ and $\gamma'$
is the unique standard parabolic Coxeter element with this property.
\begin{lemma}  If $F$ is a face of $\Delta (\gamma)$
and $F'$ is the face of $F$ obtained by removing a set of negative simple
roots from the vertex set of $F$, then the sets of right vertices of $F$
and $F'$ coincide.
\label{lem:neg-simple}
\end{lemma}
\begin{proof}
It suffices to consider the case where $F'$ is obtained from $F$
by removing a single negative simple root. Let $V =
\{\tau_1,\dots,\tau_k\}$ be the ordered vertex set of $F$. By
definition of the order (\ref{total}) there exist integers $i$ and
$j$ with $0 \le i <j \le k$ such that
\[\{\tau_1,\dots,\tau_i\} \subseteq -\Pi_2, \ \ \ \
\{\tau_{i+1},\dots,\tau_j\} \subseteq \Phi^+ \ \ \mbox{ and } \ \
\{\tau_{j+1},\dots,\tau_k\} \subseteq -\Pi_1.\]
In view of Lemma \ref{lem:easyleft} it suffices to show that for $i+1 \le
p \le j$, the type of $\tau_p$ is unchanged if a negative simple root is
removed from $V$.  Since the type of $\tau_p$ is determined
by $\{\tau_p, \tau_{p+1},\dots,\tau_k \}$ we need only consider the
removal of a negative simple root in the set $\{\tau_{j+1},\dots,\tau_k\}$.
Suppose that $\tau_q$ is such a negative simple root and that the roots
$u_p (F)(\tau_p)$ and $u'_p (F) (\tau_p)$ have different signs, where
\[u_p (F) = R(\tau_k) \cdots R(\tau_p) \ \ \mbox{ and } \ \
u'_p (F) = R(\tau_k) \cdots R(\tau_{q+1})R(\tau_{q-1}) \cdots R(\tau_p).\]
Since $\Pi_1$ is an orthonormal set we have $u'_p(F) = R(\tau_q) u_p(F)$
whence, since $\tau_q$ is a negative simple root, we conclude that
$u'_p (F)(\tau_p) = \pm \tau_q$. However this forces a linear dependence
on the set $\{\tau_p,\dots,\tau_k\}$, giving a contradiction.
\end{proof}
\begin{theorem}
The map $\phi$ is a bijection from the set of facets of $\Delta (\gamma)$
to $\nc (\gamma)$.
\label{thm:bij}
\end{theorem}
\begin{proof}
To prove surjectivity of $\phi$ let $w \preceq \gamma$ and choose the
unique standard parabolic Coxeter element $\gamma'$ of the form
(\ref{gamma'}) with respect to which $w \preceq \gamma'$ is
non-peripheral. Let $\phi'$ be the map of Definition \ref{def:phi}
corresponding to $\gamma'$. By Lemma \ref{lem:surj} we can find a facet
$F'$ of $X(\gamma')$ such that $\phi' (F') = w$. Extend $F'$ to a facet
$F$ of $\Delta (\gamma)$ by adding the negative simple roots not present
in (\ref{gamma'}) and note that $\phi(F) = w$ by Lemma \ref{lem:neg-simple}.

To prove injectivity of $\phi$ let $w \preceq \gamma$. We need to show
that there is at most one facet $F$ of $\Delta (\gamma)$ with $\phi(F) = w$.
Let $\{\tau_1,\dots,\tau_k\}$ be the ordered set of positive vertices of
such a facet $F$, forming the face $\tilde{F}$ of $F$. Denote by
$\gamma'$ and $\tilde{\gamma}$ the unique standard parabolic Coxeter
elements with respect to which $w \preceq \gamma'$ and $R(\tau_k) \cdots
R(\tau_1) \preceq \tilde{\gamma}$ are non-peripheral, respectively.
By Lemma \ref{lem:neg-simple} we have $\tilde{\phi} (\tilde{F}) =
w$, where $\tilde{\phi}$ is the map of Definition \ref{def:phi}
corresponding to $\tilde{\gamma}$. Proposition \ref{prop:per} implies that
$w$ is non-peripheral with respect to $\tilde{\gamma}$ and hence we must
have $\tilde{\gamma} = \gamma'$. Thus the set of negative
vertices of $F$ is equal to the negative of the set of simple roots not
appearing in the expression (\ref{gamma'}) for $\gamma'$ and, as a result,
this set of negative vertices is uniquely determined by $w$.
Finally, since $\tilde{\phi} (\tilde{F})
= w$, $\tilde{F}$ is uniquely determined by $w$ by Corollary \ref{cor:inj}.
\end{proof}
\begin{corollary}
The number of facets of $\Delta (\gamma)$ is equal to the number of
elements of $\nc (\gamma)$.
\qed
\label{cor:bij}
\end{corollary}
\begin{corollary}
The map $\phi$ restricts to a bijection from the set of facets of $X
(\gamma)$ to the set of non-peripheral elements of $\nc(\gamma)$.
\label{cor:non-per}
\end{corollary}
\begin{proof}
Combine Theorem \ref{thm:bij} with Proposition \ref{prop:per}.
\end{proof}

In particular, the number of facets of $X(\gamma)$, and hence of $\Delta_+
(\Phi)$, is equal to the number of non-peripheral elements of $\nc(\gamma)$.
This fact was found independently by N.~Reading \cite[Corollary 9.2]{Re}.
The number of facets of $\Delta_+ (\Phi)$ (positive clusters) is given by a
product formula similar to (\ref{prod}); see \cite[(3.8)]{FZ2}.

\section{Shellings and $h$-vectors}
\label{shell}

In this section we describe an explicit family of shellings of $\Delta
(\gamma)$ and use it to prove Theorem \ref{thm0}. We consider the
reverse of the lexicographic ordering for the various complexes under
consideration instead of the lexicographic ordering itself only because
this makes some of the statements technically easier to prove.
\begin{theorem}
The reverse of the lexicographic ordering on the facets of $X(w)$ is a
shelling of $X(w)$ for any $w \in \nc (\gamma)$.
\label{thm:shell+}
\end{theorem}
\begin{proof}
Let $F$ and $F'$ be two facets of $X(w)$ with vertex sets $V$ and $V'$,
respectively, such that $F'$ succeeds $F$ in the lexicographic order.
We need to show that there exists $\rho \in V \sm \, V'$ such that $V \sm
\, \{\rho\}$ is contained in the vertex set of a facet of $X(w)$ which
succeeds $F$ in the lexicographic order. We proceed by induction on the
length $k$ of $w$. The statement is trivial for $k \le 1$, so suppose
$k \ge 2$. Let $\tau$ and $\tau'$ be the smallest elements of $V$ and
$V'$, respectively, so that $\tau \le \tau'$. If $\tau < \tau'$ then
$\tau$ is not a vertex of the last facet of $X(w)$ in the lexicographic
order and hence $w(\tau)$ is a positive root by Corollary \ref{cor:last3}.
This means that $\tau$ is a left vertex in $F$. Clearly $\tau \notin F'$
and the result follows in this case from Proposition \ref{prop:type2} (i)
with $\rho = \tau$.
Suppose now that $\tau = \tau'$. Then $V \sm \, \{\tau\}$ and $V' \sm \,
\{\tau\}$ are the vertex sets of facets $G$ and $G'$, respectively, of
$X(w R(\tau))$ such that $G'$ succeeds $G$ in the lexicographic order.
By induction $G$ precedes a facet of $X(w R(\tau))$ with vertex
set of the form $(V \sm \, \{\tau, \rho\}) \cup \{\rho'\}$ for some $\rho
\in V \sm \, V'$, so that necessarily $\rho < \rho'$. It follows that
$(V \sm \, \{\rho\}) \cup \{\rho'\}$ is the vertex set of a facet of
$X(w)$ which succeeds $F$. This completes the induction.
\end{proof}

Let $n(F)$ denote the number of vertices of a face $F$ of $\Delta (\gamma)$
and $F_+$ (respectively, $F_-$) denote the face of $F$ whose vertices are
the positive (respectively, negative) vertices of $F$. Define the
partial order $\unlhd$ on the set of facets of $\Delta (\gamma)$ as follows.
For two such facets $F$ and $F'$ we have $F' \lhd F$ if and only if either
$n (F'_+) > n (F_+)$ or $n(F'_+) = n(F_+)$ and $F_+$ precedes $F'_+$
in the lexicographic order.
\begin{lemma}
Let $V$ be the vertex set of a facet $F$ of $\Delta (\gamma)$, let $\tau \in
V$ and let $\tau'$ be the unique vertex of $\Delta (\gamma)$ other than
$\tau$ such that $(V \sm \, \{\tau\}) \cup \{\tau'\}$ is the vertex set
of a facet of $\Delta (\gamma)$.
\begin{enumerate}
\itemsep=0pt
\item[(i)] If $\tau$ is a negative root then $\tau'$ is a positive root.

\item[(ii)] Let $\tau$ be positive. Then $\tau$ is a left vertex in
$F$ if and only if $\tau'$ is positive and $\tau < \tau'$.
\end{enumerate}
\label{lem:restriction}
\end{lemma}
\begin{proof}
(i) Observe that if $\tau$ is a negative root then the wall of $F$ opposite
$\tau$ contains all other negative simple roots.

\smallskip
\noindent
(ii) Let $V_+$ denote the set of positive vertices of $F$ and let $w \preceq
\gamma$ be such that $F_+$ is a facet of $X(w)$. By Lemma
\ref{lem:neg-simple} we have that $\tau$ is a left vertex in $F$ if and
only if it is a left vertex in $F_+$. By Proposition \ref{prop:type2} (i)
this happens if and only if there exists a positive root $\sigma > \tau$
such that $(V_+ \sm \, \{\tau\}) \cup \{\sigma\}$ is the vertex set of a
facet of $X(w)$. To complete the proof observe that for $\sigma$ positive,
$(V_+ \sm \, \{\tau\})
\cup \{\sigma\}$ is the vertex set of a facet of $X(w)$ if and only if
$(V \sm \, \{\tau\}) \cup \{\sigma\}$ is the vertex set of a facet of
$\Delta (\gamma)$.
\end{proof}

\begin{theorem}
\begin{enumerate}
\itemsep=0pt
\item[(i)] Any linear extension of the partial order $\unlhd$ is a
shelling order for $\Delta (\gamma)$.

\item[(ii)] The restriction set of a facet $F$ of $\Delta (\gamma)$ with
respect to this shelling is equal to the set of left vertices in $F$.
\end{enumerate}
\label{thm:shell}
\end{theorem}
\begin{proof}
(i) Let $F$ and $F'$ be two facets of $\Delta (\gamma)$ with vertex sets $V$
and $V'$, respectively, such that $F' \lhd F$. It suffices to show that there
exists $\rho \in V \sm \, V'$ such that $V \sm \, \{\rho\}$ is contained in
the vertex set of a facet of $\Delta (\gamma)$ which precedes $F$ in the
order $\unlhd$. If $F'_- = F_-$ then the statement follows from Theorem
\ref{thm:shell+}. Otherwise one can choose as $\rho$ any negative root in
$V \sm \, V'$ since then, by Lemma \ref{lem:restriction} (i), $V \sm \,
\{\rho\}$ is contained in the vertex set of a facet of $\Delta (\gamma)$
having one more positive vertex than $F$.

\smallskip
\noindent
(ii) This is an immediate consequence of the definition of the restriction
set and of Lemmas \ref{lem:easyleft} and \ref{lem:restriction}.
\end{proof}

%%total ordering
%
%%\[ -\alpha_1,\dots,-\alpha_n, \, \rho_1, \rho_2,\dots,\rho_{nh/2} \]
%
%%of the vertex set of $GA(W)$.
%
%%\begin{theorem}
%%The reverse of the lexicographic ordering on the facets of $GA(W)$ is
%%a shelling.
%%\label{thm:shell}
%%\end{theorem}
%
%%\begin{proof}
%%For $0 \le i \le n$ let $X_i (\gamma)$ denote the induced subcomplex of
%%$GA(W)$ on the vertex set $P \cup \{-\alpha_n,\dots,-\alpha_{n-i+1}\}$,
%%so that $X_0 (\gamma) = X (\gamma)$ and $X_n (\gamma) = GA(W)$. We will
%%show that the reverse of the lexicographic ordering on the facets of $X_i
%%(\gamma)$ is a shelling for all $i$. We proceed by induction on $i$ and
%%the rank of $W$. For $i=0$ the statement reduces to Theorem
%%\ref{thm:shell+}. For $i \ge 1$ let $j = n-i+1$ and observe that $X_i
%%(\gamma) \sm \, \{-\alpha_j\} = X_{i-1} (\gamma)$ and $X_i (\gamma) /
%%\{-\alpha_j\} = X_{i-1} (\gamma')$ where $\gamma' = R(\alpha_1) \cdots
%%R(\alpha_{j-1}) R(\alpha_{j+1}) \cdots R(\alpha_n)$, so the statement
%%follows from the induction hypotheses and Lemma ??.
%%\end{proof}

%
\begin{corollary}
\begin{enumerate}
\itemsep=0pt
\item[(i)] The entry $h_i (\Delta (\gamma))$ of the $h$-vector of $\Delta
(\gamma)$ is equal to the number of elements of $\nc(\gamma)$ of rank $i$.

\item[(ii)] The entry $h_i (X(\gamma))$ of the $h$-vector of $X(\gamma)$
is equal to the number of non-peripheral elements of $\nc (\gamma)$ of
rank $n-i$.
\end{enumerate}
\label{cor:h-vector}
\end{corollary}
\begin{proof}
(i) It follows from Theorem \ref{thm:shell} and (\ref{h-shell})
that $h_i (\Delta (\gamma))$ is equal to the number of facets of $\Delta
(\gamma)$ with $i$ left vertices or, equivalently, to the number of
facets of $\Delta (\gamma)$ with
$n-i$ right vertices. By Theorem \ref{thm:bij} this is equal to the
number of elements of $\nc(\gamma)$ of rank $n-i$, which is equivalent
to the statement we want to prove by the self-duality of $\nc (\gamma)$.

\noindent
(ii) The facets of $X(\gamma)$ come first in the shelling of Theorem
\ref{thm:shell} and hence $h_i (X(\gamma))$ is equal to the number of
facets of $X(\gamma)$ with $i$ left vertices or, equivalently, to the
number of facets of $X(\gamma)$ with $n-i$ right vertices. This is
equal to the number of non-peripheral elements of $\nc (\gamma)$ of
rank $n-i$ by Corollary \ref{cor:non-per}.
\end{proof}

Part (ii) of the previous corollary gives a combinatorial
interpretation of the entries of the $h$-vector of $X(\gamma)$,
and hence of $\Delta_+ (\Phi)$. Several different combinatorial
interpretations to these numbers appear in \cite{AT}. The
following statement, which is a direct consequence of Proposition
\ref{prop:per-char} and part (ii) of Corollary \ref{cor:h-vector},
provides a similar interpretation.
\begin{corollary}
The entry $h_i (X(\gamma))$ of the $h$-vector of $X(\gamma)$ is equal
to the number of elements of $\nc (\gamma)$ of rank $i$ which are not
preceded by any of the reflections corresponding to the last $n$
positive roots.
\qed
\label{cor:h-vector+}
\end{corollary}
\begin{remark}
{\rm  The non-peripheral elements of $\nc (\gamma)$ of rank one
are the reflections in $W$ which do not lie in any proper standard
parabolic subgroup of $W$. The number $f_W$ of such reflections
was studied by F. Chapoton~\cite{Ch2} and can be expressed in
terms of the exponents and Coxeter number of $W$ \cite[Proposition
1.1]{Ch2}. Hence in the special case $i=n-1$, part (ii) of
Corollary \ref{cor:h-vector} provides a conceptual explanation of
the equality $h_{n-1} (\Delta_+ (\Phi)) = f_W$, which can be
checked case by case on the basis of the data provided in
\cite[Section 6]{AT} and \cite[Section 1]{Ch2}.} \label{rem:ch2}
\end{remark}

\bigskip
\noindent \emph{Proof of Theorem \ref{thm0}.} Let $\Phi$ and $W$
be as in the statement of the theorem. If $\Phi = \Phi_1 \times
\Phi_2 \times \cdots \times \Phi_m$ is the decomposition of $\Phi$
into irreducible components and $W = W_1 \times W_2 \times \cdots
\times W_m$ is the corresponding decomposition of $W$ then $\Delta
(\Phi)$ is the simplicial join of the $\Delta (\Phi_i)$ and
$\nc_W$ is the direct product of the $\nc_{W_i}$ and hence the
$h$-polynomial of $\Delta (\Phi)$ (the polynomial with
coefficients the entries of the $h$-vector) and the rank
generating polynomial of $\nc_W$ are multiplicative with respect
to these decompositions. Therefore we may assume that $\Phi$
(equivalently, $W$) is irreducible. Since the combinatorial
structure of $\Delta (\Phi)$ is unaffected by rescaling of roots,
we may assume further that the elements of $\Phi$ are of unit
length. Under these assumptions the result follows from Corollary
\ref{cor:h-vector} (i) and Theorem \ref{thm:BW} (i). \qed

\vspace{0.2 in} \noindent \emph{Acknowledgements}. The first
author thanks the American Institute of Mathematics, Palo Alto,
and the Mittag-Leffler Institute, Dj\"ursholm, Sweden, for their
hospitality and financial support. The second author thanks the
American Institute of Mathematics, Palo Alto, the Centre de
Recerca Matem\`{a}tica, Barcelona, and the Centre de
Math\'{e}matiques et Informatique, Marseille, for their
hospitality and financial support.  The third author would like to
thank the American Institute of Mathematics and the National Science
Foundation for their hospitality and support.

\end{document}